\documentclass[12pt]{article}
\usepackage{graphics,graphicx}
\usepackage{amssymb,amsmath,amsopn,amsfonts}
\usepackage[dvips]{color}
\usepackage{colordvi,multicol}
\usepackage{epsf}
\newfam\msbfam
\font\tenmsb=msbm10 \textfont\msbfam=\tenmsb \font\sevenmsb=msbm7
\scriptfont\msbfam=\sevenmsb \font\fivemsb=msbm5
\scriptscriptfont\msbfam=\fivemsb

\def\th#1{\vspace{1mm}\noindent{\bf #1}\quad}

\def\proof{\vspace{1mm}\noindent{\it Proof}\quad}

\numberwithin{equation}{section}

\def\bc{\begin{center}}
\def\ec{\end{center}}
\def\no{\noindent}
\def\hang{\hangindent\parindent}
\def\textindent#1{\indent\llap{\qquad #1\ \ \enspace}\ignorespaces}
\def\ref{\par\hang\textindent}

\begin{document}
\title{ {\bf Piecewise linear approximation for  the dynamical $\Phi^4_3$ model
\thanks{Research supported in part  by NSFC (No.11671035, No.11401019)}\\} }
\author{  {\bf Rongchan Zhu}$^{\mbox{a,c}}$, {\bf Xiangchan Zhu}$^{\mbox{b,c},}$\thanks{Corresponding author}
\date{}
\thanks{E-mail address:
zhurongchan@126.com(R. C. Zhu), zhuxiangchan@126.com(X. C. Zhu)}\\\\
$^{\mbox{a}}$Department of Mathematics, Beijing Institute of Technology, Beijing 100081,  China\\
$^{\mbox{b}}$School of Science, Beijing Jiaotong University, Beijing 100044, China\\
$^{\mbox{c}}$Department of Mathematics, University of Bielefeld, D-33615 Bielefeld, Germany}

\maketitle

\noindent {\bf Abstract} We construct a piecewise linear approximation for  the dynamical $\Phi_3^4$ model on $\mathbb{T}^3$ by the theory of regularity structures in [Hai14].
For the dynamical $\Phi^4_3$ model it is proved in [Hai14] that a renormalisation has to be performed in order to define the nonlinear term. Compared to the results in [Hai14] we consider  piecewise linear approximations to  space-time white noise and prove that the solutions to the approximating equations converge to the solution to the dynamical $\Phi^4_3$ model. The renormalisation in this case corresponds to adding the solution multiplied by a function depending on $t$ in the approximating equation.

\vspace{1mm}
\no{\footnotesize{\bf 2000 Mathematics Subject Classification AMS}:\hspace{2mm} 60H15, 82C28}
 \vspace{2mm}

\no{\footnotesize{\bf Keywords}:\hspace{2mm}   regularity structure,  space-time white noise, SPDEs, renormalisation}

\section{Introduction}

In this paper we construct a piecewise linear approximation of the dynamical $\Phi^4_3$ model driven by space-time white noise on $\mathbb{T}^3$:
\begin{equation}d\Phi =\Delta \Phi dt-\Phi^3 dt +  dW(t).\end{equation}
This can be considered as a Wong-Zakai approximation (c.f. [WZ65a,WZ65b]) for additive noise in this case.
Here $W$ is a two-sided  cylindrical Wiener process on $L^2(\mathbb{T}^3)$. Formally, if we define $\xi$ by $\int\phi(t,x)\xi(dt,$ $dx)=\int\langle\phi,dW(s)\rangle$ for $\phi\in L^2(\mathbb{R}\times \mathbb{T}^3)$, then $\xi$ is  periodic space-time white noise.  This model is also known as the stochastic quantisation of
Euclidean quantum field theory ( see [GJ87] and the reference therein). It is also considered as a universal
model for phase coexistence near the critical point (see [GLP99]).
In two spatial dimensions, this problem was previously treated in [AR91] and [DD03].
In three spatial dimensions  this  equation (1.1)  is ill-posed and  the main difficulty in this case is that $W$ and hence $\Phi$ are so singular that the non-linear term is not well-defined in the classical sense. It was a long-standing open problem to give a meaning to this equation in the three dimensional case.

A breakthrough result was achieved recently by Martin Hairer in [Hai14], where he introduced a theory of regularity structures and gave a meaning to this equation (1.1) successfully. Also by using the paracontrolled distribution method proposed by Gubinelli, Imkeller and  Perkowski in [GIP13]  existence and uniqueness of local solutions to (1.1) has been obtained in [CC13]. Recently, these two approaches have been successful in giving a meaning to a lot of ill-posed stochastic PDEs like the Kardar-Parisi-Zhang (KPZ) equation ([KPZ86], [BG97], [Hai13]), the dynamical $\Phi_3^4$ model ([Hai14], [CC13]), the Navier-Stokes equation driven by space-time white noise ([ZZ14], [ZZ15]), the dynamical sine-Gordon equation ([HS14]) and so on (see [HP14] for more interesting examples).  From a “philosophical” perspective, the theory of regularity structures and  the paracontrolled distribution  are inspired by the theory
of controlled rough paths (see [Lyo98], [Gub04]). The main difference is that the regularity structure theory considers the problem locally, while the paracontrolled distribution method is a global approach  using Fourier analysis.

An interesting question for SDE, especially the dynamical $\Phi^4_3$ model is as follows: Given a sequence $W_\varepsilon$ of regularization of the  noise $W$ (for example convolutions with a mollifier), can we obtain a non-trivial  solution associated with $W$ by taking the limit of $\Phi_\varepsilon$ as $\varepsilon$ goes to $0$, where $\Phi_\varepsilon$ is the solution associated to $W_\varepsilon$. In the finite dimensional case a series of classical results has been obtained by  Wong and Zakai [WZ65a,WZ65b].  However, the answer to this question for  the dynamical $\Phi^4_3$ model is no (see [HRW12]). Indeed, we have to consider the following modified equation
\begin{equation}\partial_t\Phi_\varepsilon =\Delta\Phi_\varepsilon+C_\varepsilon\Phi_\varepsilon-\Phi_\varepsilon^3+\xi_\varepsilon.\end{equation}
In [Hai14] Martin Hairer considered an $\varepsilon$-approximation $\xi_\varepsilon$ to space-time white noise. Here  $\xi_\varepsilon$ is given by  convolution with a mollifier, i.e. $\xi_\varepsilon:=\rho_\varepsilon*\xi$, where the convolution means that we view $\xi$ as a distribution on $\mathbb{R}^4$ and do convolution on $\mathbb{R}^4$. $\rho_\varepsilon(t,x)$ is a compactly
supported smooth mollifier  that is scaled by $\varepsilon$ in the spatial directions and by $\varepsilon^2$ in the time
direction, i.e. $\rho_\varepsilon(t,x)=\varepsilon^{-5}\rho(\varepsilon^{-2}t,\varepsilon^{-1}x)$ for some smooth, compactly supported function $\rho$. Let $\Phi_\varepsilon$ denote the solution to (1.2).
It is proved in [Hai14] that there exist choices of constants $C_\varepsilon$ diverging as $\varepsilon\rightarrow0$, as well as a process $\Phi$ such that $\Phi_\varepsilon\rightarrow \Phi$ in probability. Furthermore, while the constants $C_\varepsilon$  do
depend crucially on the choice of the mollifiers $\rho_\varepsilon$, the limit $\Phi$ does not depend on them. Also in [CC13]  purely spatial regularization has been considered and a  similar result has been obtained.

In this paper we consider another approximation given by   piecewise linear approximations combined with  convolution
with a mollifier. First, we convolute with a mollifier: $W_\varepsilon(t)=\int_0^t\xi_\varepsilon(s)ds, t\in\mathbb{R}$, with $\xi_\varepsilon$ given as above and then
consider  piecewise linear approximations: for $t\in[k\vartheta, (k+1)\vartheta)$, $k\in\mathbb{Z}$, $\vartheta>0$ $$W_{\varepsilon,\vartheta}(t)=W_\varepsilon(k\vartheta)+\frac{t-k\vartheta}{\vartheta}(W_\varepsilon((k+1)\vartheta)-W_\varepsilon(k\vartheta)),$$
and $\xi_{\varepsilon,\vartheta}(t)=\partial_tW_{\varepsilon,\vartheta}=\frac{1}{\vartheta}\int_{k\vartheta}^{(k+1)\vartheta}\xi_{\varepsilon}(u)du$ for $t\in[k\vartheta, (k+1)\vartheta)$, $k\in\mathbb{Z}$, which is our regularised noise.

 This approximation is  the
celebrated Wong-Zakai approximation of the solution and is related to a classical problem:  approximating  solutions in terms of a simpler
model, where the stochastic integral is changed into a ”deterministic” one and replacing the
noise by  its piecewise linear interpolation on a time grid.  For  finite-dimensional diffusion processes, this kind of approximation is well-known
(see, e.g. [T96, LQZ02]  and the references therein). There is
a substantial number of publications devoted to Wong-Zakai approximations of infinite dimensional
stochastic equations (see [N04, CM11] and the references therein).

In this paper we  use the theory of regularity structures to study this approximation for the dynamical $\Phi_3^4$ model. The key idea of the theory of regularity structures is as follows:  we perform an abstract Talyor expansion on  both sides of the equation. Originally, Talyor expansions are only for functions. Here  the right objects, e.g. a regularity structure whose elements could
possibly take the place of Taylor polynomials, can be constructed. Given a noise $\xi$, the regularity structure can  be endowed with a model $\iota\xi$, which is a concrete way of associating every element in  the abstract regularity structure to the actual Taylor polynomial at every point. Multiplication, differentiation, the state space of solutions,  and the convolution with singular kernels can be defined on this regularity structure, which is the major difficulty when trying to give a meaning to such singular stochastic partial differential equations as above. On the regularity   structure, a fixed point argument can be applied to obtain  local existence and uniqueness of the solutions $\bar{\Phi}$ to the equation lifted onto the regularity structure. Furthermore, we can go back to "the real world" with the help of another central tool of the theory, namely the reconstruction operator $\mathcal{R}$. If $\xi$ is a smooth process, $\Phi=\mathcal{R}\bar{\Phi}$ coincides with the classic solution to the equation. Now we have the following maps
$$\xi\mapsto\iota\xi\mapsto\bar{\Phi}\mapsto\mathcal{R}\bar{\Phi}.$$
The last two maps are continuous with respect to suitable topologies, while   the above sequence $\iota\xi_{\varepsilon}$ of
canonical models fails to converge with a smooth approximation $\xi_\varepsilon$ to the noise $\xi$. It may, however,  still be possible to renormalise the model $\iota\xi_{\varepsilon}$ to some converging
model $\hat{\iota}\xi_{\varepsilon}$, which in turn can be related to a specific renormalised equation (1.2).

In this paper for the approximating  sequence $\xi_{\varepsilon,\vartheta}$ we build the associated model $\iota\xi_{\varepsilon,\vartheta}$, which also need to be renormalised  into some converging
model $\hat{\iota}\xi_{\varepsilon,\vartheta}$. This in turn can be related to the following renormalised equation (1.3):
\begin{equation}\partial_t\Phi_{\varepsilon,\vartheta}(t) =\Delta\Phi_{\varepsilon,\vartheta}(t)+C^{(\varepsilon,\vartheta)}(t)\Phi_{\varepsilon,\vartheta}(t)-\Phi_{\varepsilon,\vartheta}^3(t)+\xi_{\varepsilon,\vartheta}(t).\end{equation}
 Here $C^{(\varepsilon,\vartheta)}$ are functions depending only on time $t$.

With these notations at hand, the main result of this article is as follows:

 \vskip.10in
 \th{Theorem 1.1} Let $\xi_{\varepsilon,\vartheta}$ be defined as above.   Denote by $\Phi_{\varepsilon,\vartheta}$ the solution to (1.3). Suppose that $\rho(t,x)=\rho_1(t)\rho_2(x)$ for smooth functions $\rho_1, \rho_2$. Then there exist choices of functions $C^{(\varepsilon,\vartheta)}$ diverging as $\varepsilon, \vartheta\rightarrow0$ such that $\Phi_{\varepsilon,\vartheta}\rightarrow \Phi$ in probability locally in time. Here  $\Phi$ is  the solution to the dynamical $\Phi^4_3$ model obtained in [Hai14].

\vskip.10in

 \th{Remark 1.2} (i)  We can also first do  purely spatial regularization corresponding to $\rho(t,x)=\delta(t)\rho_2(x)$ for the Dirac distribution $\delta$ and then do  piecewise linear approximation. In this case the results in Theorem 1.1 still holds (see Remark 3.8). In fact, the only difference is the proof of Theorem 3.7.

(ii) The function $C^{(\varepsilon,\vartheta)}$ in (1.3) is given as follows:
$$C^{(\varepsilon,\vartheta)}:=3C_{1}^{(\varepsilon,\vartheta)}-9C_{2}^{(\varepsilon,\vartheta)},$$ where $C_1^{(\varepsilon,\vartheta)}$ and $C_2^{(\varepsilon,\vartheta)}$ are defined in (A.4) and (A.8).   Moreover, $$C_1^{(\varepsilon,\vartheta)}\backsimeq
\frac{1}{\varepsilon},\quad C_{2}^{(\varepsilon,\vartheta)}\backsimeq-\log\varepsilon.$$
Here we emphasize that each function $C_i^{(\varepsilon,\vartheta)}(t), i=1,2,$ cannot be separated as a diverging constant  and a converging function. In fact, for the case that $\rho(t,x)=\delta(t)\rho_2(x)$, by a straightforward calculation we know that $|C_1^{(\varepsilon,\vartheta)}(t_1)-C_1^{(\varepsilon,\vartheta)}(t_2)|\simeq\frac{1}{\varepsilon}$ for $t_1\neq t_2$.

\vskip.10in
As mentioned in Remark 1.2 (ii), in our case it is required in (1.3) to minus $\Phi_{\varepsilon,\vartheta}$ multiplied by a function $C^{(\varepsilon,\vartheta)}$ depending on $t$ such that the associated solutions $\Phi_{\varepsilon,\vartheta}$ converge to the solution to the dynamical $\Phi^4_3$ model  as $\varepsilon,\vartheta\rightarrow0$, which is the main difference from the result in [Hai14]. We introduce a new symbol $\mathbf{C}$ in the regularity structure to represent $C^{(\varepsilon,\vartheta)}(t)$ and  define a bigger regularity structure $\mathfrak{T}^1$  including $\mathbf{C}$ as well as the original regularity structure $\mathfrak{T}_F$  associated with the $\Phi^4_3$ model, which helps us  to construct  a suitable renormalised model corresponding to (1.3) for $\mathfrak{T}_F$ (see Remark 3.4).

We would also like to emphasize that the proof  in this paper is
not restricted to the specific equation (1.1). Similar arguments would yield similar
results for  the models  that can be treated with the methods developed in [Hai14].

In Section 2 we present a  summary of
some notions of the theory of regularity structures. In Section 3 we construct the renormalised model and prove the main results. The convergence of the renormalised model is proved in Section 4. The Appendix contains the  proof of Theorem 3.7.

\section{Regularity structures}
 In this section we recall some preliminaries for the  theory of regularity structures from [Hai14].
\vskip.10in
\th{Definition 2.1} A regularity structure $\mathfrak{T}=(A,T,G)$ consists of the following elements:

(i) An index set $A\subset \mathbb{R}$ such that $0\in A$, $A$ is bounded from below and  locally finite.

(ii) A model space $T$, which is a graded vector space $T=\oplus_{\alpha\in A}T_\alpha$, with each $T_\alpha$ a Banach space. Furthermore, $T_0$ is one-dimensional and has a basis vector $\mathbf{1}$. Given $\tau\in T$ we write $\|\tau\|_\alpha$ for the norm of its component in $T_\alpha$.

(iii) A structure group $G$ of (continuous) linear operators acting on $T$ such that for every $\Gamma\in G$, every $\alpha\in A$ and every $\tau_\alpha\in T_\alpha$ one has
$$\Gamma\tau_\alpha-\tau_\alpha\in T_{<\alpha}:=\bigoplus_{\beta<\alpha}T_\beta.$$
Furthermore, $\Gamma\mathbf{1}=\mathbf{1}$ for every $\Gamma\in G$.
\vskip.10in
The canonical example is the space $\bar{T} =\bigoplus_{n\in\mathbb{N}}\bar{T}_n$ of abstract
polynomials in finitely many indeterminates, with $A = \mathbb{N}$ and $\bar{T}_n$ denoting the
space of monomials that are homogeneous of degree $n$. In this case, a natural group
of transformations $G$ acting on $\bar{T}$ is given by the group of translations.

Given a scaling $\mathfrak{s}=(\mathfrak{s}_0,\mathfrak{s}_1,...,\mathfrak{s}_{d})$ of $\mathbb{R}^{d+1}$. We call $|\mathfrak{s}|=\mathfrak{s}_0+\mathfrak{s}_1+...+\mathfrak{s}_d$ scaling dimension.  We define the  associate  metric on $\mathbb{R}^{d+1}$  by
$$\|z-z'\|_\mathfrak{s}:=d_\mathfrak{s}(z,z'):=\sum_{i=0}^{d}|z_i-z'_i|^{1/\mathfrak{s}_i}.$$
For $k=(k_0,...,k_d)$ we set $|k|_\mathfrak{s}=\sum_{i=0}^d\mathfrak{s}_ik_i$.

\subsection{Specific regularity structures}

We start with the regularity structure $\bar{T}$ given by all polynomials in $d+1$ indeterminates, let us call them $X_0,..., X_d$, which denote the time and space directions respectively. Recall $X^k=X_0^{k_0}\cdot\cdot\cdot X_d^{k_d}$ with $k$ a multi-index.
For the case of the dynamical  $\Phi^4_3$ model, $d=3$ and the scaling is $\mathfrak{s}=(2,1,1,1)$.  In the regularity structure we use the symbol  $\Xi$ to replace the driving noise $\xi$. We introduce the integration maps $\mathcal{I}$ and $\mathcal{I}_k$ for a multi-index $k$ associated with the operation of convolution with a truncation of the heat kernel $G$ and
its  derivative respectively.

We recall the following notations from [Hai14]: define a set $\mathcal{F}$ by postulating that $\{\mathbf{1},\Xi,X_j\}\subset \mathcal{F}$ and whenever $\tau,\bar{\tau}\in\mathcal{F}$, we have
$\tau\bar{\tau}\in\mathcal{F}$ and $\mathcal{I}_k(\tau)\in\mathcal{F}$; define $\mathcal{F}_+$ as the set of all elements $\tau\in\mathcal{F}$ such that either $\tau=\mathbf{1}$ or $|\tau|_\mathfrak{s}>0$ and such that, whenever $\tau$ can be written as $\tau=\tau_1\tau_2$ we have either $\tau_i=\mathbf{1}$ or $|\tau_i|_\mathfrak{s}>0$; $\mathcal{H}, \mathcal{H}_+$ denote the sets of finite linear combinations of all elements in $\mathcal{F}, \mathcal{F}_+$, respectively.
Here for each $\tau\in \mathcal{F}$ a weight $|\tau|_{\mathfrak{s}}$ is obtained by setting
$|\mathbf{1}|_\mathfrak{s}=0$,
$$|\tau\bar{\tau}|_\mathfrak{s}=|\tau|_\mathfrak{s}+|\bar{\tau}|_\mathfrak{s},$$
for any two formal expressions $\tau$ and $\bar{\tau}$ in $\mathcal{F}$ such that
$$|\Xi|_\mathfrak{s}=\alpha,\quad |X_i|_\mathfrak{s}=\mathfrak{s}_i,\quad |\mathcal{I}_k(\tau)|_\mathfrak{s}=|\tau|_\mathfrak{s}+2-|k|_\mathfrak{s},$$
 with $ -\frac{18}{7}<\alpha<-\frac{5}{2}$.

As in [Hai14] we construct the regularity structure, which contains those that are actually useful for the abstract reformulation of the equation (1.1). Define
$$\mathfrak{M}_F=\{\Xi, U^n:n\leq 3\},$$
and the sets $\mathcal{W}_0=\mathcal{U}_0=\emptyset$ and $\mathcal{W}_n, \mathcal{U}_n$ for $n>0$ recursively by
$$\mathcal{W}_n=\mathcal{W}_{n-1}\cup\bigcup_{\mathcal{Q}\in\mathfrak{M}_F}\mathcal{Q}(\mathcal{U}_{n-1},\Xi),$$
$$\mathcal{U}_n=\{X^k\}\cup\{\mathcal{I}(\tau):\tau\in\mathcal{W}_n\},$$
and
$$\mathcal{F}_F:=\bigcup_{n\geq0}(\mathcal{W}_n\cup\mathcal{U}_n).$$

Then $\mathcal{F}_F$ contains the elements required to describe both the solutions and the terms  in the equation.
We denote by $\mathcal{H}_F$ the set of all finite linear combinations of elements in $\mathcal{F}_F$.
\vskip.10in

Now we follow [Hai14] to construct the structure group $G_F$. Define a linear projection operator $P_+:\mathcal{H}\rightarrow\mathcal{H}_+$ by imposing that
$$P_+\tau=\tau,\quad \tau\in\mathcal{F}_+,\quad P_+\tau=0,\quad \tau\in\mathcal{F}\setminus\mathcal{F}_+,$$
and two linear maps $\Delta:\mathcal{H}\rightarrow\mathcal{H}\otimes \mathcal{H}_+$ and $\Delta^+:\mathcal{H}_+\rightarrow\mathcal{H}_+\otimes \mathcal{H}_+$
by
$$\Delta\mathbf{1}=\mathbf{1}\otimes \mathbf{1},\quad \Delta^+\mathbf{1}=\mathbf{1}\otimes \mathbf{1},$$
$$\Delta X_i=X_i\otimes \mathbf{1}+\mathbf{1}\otimes X_i, \quad \Delta^+ X_i=X_i\otimes \mathbf{1}+\mathbf{1}\otimes X_i,$$
$$\Delta \Xi=\Xi\otimes \mathbf{1},$$
and recursively by $$\Delta(\tau\bar{\tau})=(\Delta\tau)(\Delta\bar{\tau})$$
$$\Delta(\mathcal{I}_k\tau)=(\mathcal{I}_{k}\otimes I)\Delta\tau+\sum_{l,m}\frac{X^l}{l!}\otimes \frac{X^m}{m!}(P_+\mathcal{I}_{k+l+m}\tau),$$
$$\Delta^+(\tau\bar{\tau})=(\Delta^+\tau)(\Delta^+\bar{\tau})$$
$$\Delta^+(\mathcal{I}_k\tau)=(I\otimes \mathcal{I}_{k}\tau)+\sum_{l}(P_+\mathcal{I}_{k+l}\otimes \frac{(-X)^l}{l!})\Delta\tau.$$

 By using the theory of regularity structures (see [Hai14, Section 8])  a structure group $G_F$ of linear operators acting on  $\mathcal{H}_F$ satisfying Definition 2.1 can be defined as follows: For  $g\in \mathcal{H}_+^*$, the dual of $\mathcal{H}_+$, satisfying $g(\tau \bar{\tau})=g(\tau)g(\bar{\tau})$ for $\tau,\bar{\tau}\in\mathcal{H}_+$, define  $\Gamma_g:\mathcal{H}\rightarrow\mathcal{H}, \Gamma_g\tau=(I\otimes g)\Delta\tau$.  Following [Hai14, Theorem 8.24]  the  regularity structure associated with the dynamical $\Phi^4_3$ model can be constructed:

 Let $T=\mathcal{H}_F$ with $T_\gamma=\langle\{\tau\in\mathcal{F}_F:|\tau|_\mathfrak{s}=\gamma\}\rangle$, $A=\{|\tau|_\mathfrak{s}:\tau\in\mathcal{F}_F\}$. Then $\mathfrak{T}_F=(A,\mathcal{H}_F,G_F)$ defines a regularity structure  associated with the dynamical $\Phi^4_3$ model.

\subsection{Models}
Now that we have fixed our algebraic regularity structure $\mathfrak{T}_F=(A,\mathcal{H}_F,G_F)$, we introduce a
family of  objects which is a concrete way of associating  every $\tau\in \mathcal{H}_F$ and $x_0\in\mathbb{R}^{d+1}$ with
the actual "Taylor polynomial based at $x_0$" represented by $\tau$
 in order to allow us to describe solutions to (1.1) locally.

First we introduce some notations:
Given a smooth compactly supported test function $\varphi$ and a space-time coordinate $z=(t,x_1,...,x_d)\in\mathbb{R}^{d+1}$, we write $\varphi_z^\lambda$ as a shorthand for
$$\varphi_z^\lambda(s,y_1,...,y_d)=\lambda^{-|\mathfrak{s}|}\varphi(\frac{s-t}{\lambda^{\mathfrak{s}_0}},\frac{y_1-x_1}{\lambda^{\mathfrak{s}_1}},...,
\frac{y_d-x_d}{\lambda^{\mathfrak{s}_d}}).$$

Let $\mathcal{B}_\alpha$ denote the set of all smooth test functions $\varphi:\mathbb{R}^{d+1}\mapsto\mathbb{R}$ that are supported in the centred  ball of radius $1$ and such that their derivatives of order up to $1+|\alpha|$ are uniformly bounded by $1$. We also denote by $\mathcal{S}'$ the space of all distributions on $\mathbb{R}^{d+1}$ and denote by $L(E,F)$ the set of all continuous linear maps between the topological vector spaces $E$ and $F$. With these notations at hand we give the definition of a model:
\vskip.10in
\th{Definition 2.2} Given a regularity structure $\mathfrak{T}=(A,T,G)$, a model for $\mathfrak{T}$ consists of maps
$$\mathbb{R}^{d+1}\ni z\mapsto\Pi_z\in L(T,\mathcal{S}'), \quad \mathbb{R}^{d+1}\times\mathbb{R}^{d+1}\ni(z,z')\mapsto\Gamma_{zz'}\in G,$$
satisfying the algebraic compatibility conditions
$$\Pi_z\Gamma_{zz'}=\Pi_{z'}, \quad \Gamma_{zz'}\circ\Gamma_{z'z''}=\Gamma_{zz''},$$
as well as the analytical bounds
$$|\Pi_z\tau(\varphi_z^\lambda)|\lesssim\lambda^\alpha\|\tau\|_\alpha, \quad \|\Gamma_{zz'}\tau\|_\beta\lesssim\|z-z'\|_\mathfrak{s}^{\alpha-\beta}\|\tau\|_\alpha.$$
Here, the bounds are imposed uniformly over all $\tau\in T_\alpha$, all $\beta<\alpha\in A$ with $\alpha<\gamma$, $\gamma>0$, and all test functions $\varphi\in \mathcal{B}_r$ with $r=\inf A$.
They are imposed locally uniformly in $z$ and $z'$.
\vskip.10in
Then for every compact set $\mathfrak{R}\subset \mathbb{R}^{d+1}$ and any two models $Z=(\Pi,\Gamma)$ and $\bar{Z}=(\bar{\Pi}, \bar{\Gamma})$ we define
$$|||Z;\bar{Z}|||_{\gamma;\mathfrak{R}}:=\sup_{z\in\mathfrak{R}}[\sup_{\varphi,\lambda,\alpha,\tau}\lambda^{-\alpha}|(\Pi_z\tau-\bar{\Pi}_z\tau)
(\varphi_z^\lambda)|+\sup_{\|z-z'\|_\mathfrak{s}\leq1}\sup_{\alpha,\beta,\tau}\|z-z'\|_\mathfrak{s}^{\beta-\alpha}
\|\Gamma_{zz'}\tau-\bar{\Gamma}_{zz'}\tau\|_\beta],$$
where the suprema are taken over the same sets as in Definition 2.2, but with $\|\tau\|_\alpha=1$. This gives a natural topology for the space  of
all models for a given regularity structure.
\vskip.10in

 To describe the models for the regularity $\mathfrak{T}_F$ we are interested in, we fix a kernel $K:\mathbb{R}^{4}\rightarrow\mathbb{R}$ with the following properties:

 (i) $K=\sum_{n\geq0}K_n$, where each  $K_n:\mathbb{R}^{4}\rightarrow\mathbb{R}$ is smooth and compactly supported in a ball of radius $2^{-n}$ around the origin. Furthermore, we assume that for every multi-index $k$, one has a constant $C$ such that
$$\sup_x|D^kK_n(x)|\leq C2^{n(2+|k|_\mathfrak{s})},$$
holds uniformly in $n$. Finally, we suppose that $\int K_n(x)P(x)dx=0$ for every polynomial $P$ of degree at most $r$ for some sufficiently large value of $r$.

(ii) $K(t,x)=0$ for $t\leq 0$ and $K(t,-x)=K(t,x)$.

(iii) For $(t,x)$ with $|x|^2+t<1/2$ and $t>0$, $K(t,x)=\frac{1}{|4\pi t|^{3/2}}e^{-\frac{|x|^2}{4t}}$, and $K$ is smooth on  $\{|x|^2+t\geq 1/4\}$.
\vskip.10in

The  kernel $K$ satisfying these properties can be obtained from the heat kernel $G$ as in [Hai14, Lemma 5.5].
\vskip.10in

\th{Definition 2.3} A model $(\Pi,\Gamma)$ for $\mathfrak{T}_F$ is admissible if it satisfies $(\Pi_xX^k)(y)=(y-x)^k$ as well as
\begin{equation}(\Pi_x\mathcal{I}\tau)(y)=\int K(y-z)(\Pi_x\tau)(z)dz+\sum_{l}\frac{(y-x)^l}{l!}f_x(P_+\mathcal{I}_{l}\tau),\end{equation}
for $\tau\in\mathcal{H}_F$ with $\mathcal{I}(\tau)\in\mathcal{H}_F$.
Here $f_x(\mathcal{I}_{l}\tau)$ are defined by
\begin{equation}f_x(\mathcal{I}_{l}\tau)=-\int D_1^lK(x-z)(\Pi_x\tau)(z)dz.\end{equation}
Furthermore, we impose $f_x(X_i)=-x_i$, $f_x(\tau\bar{\tau})=f_x(\tau)f_x(\bar{\tau})$ and extend this to all of $\mathcal{H}_+$ by linearity. $\Gamma$ is given by
\begin{equation}\Gamma_{xy}=(\Gamma_{f_x})^{-1}\circ \Gamma_{f_y},\end{equation}
where $\Gamma_{f_x}\tau:=(I\otimes {f_x})\Delta \tau$ for $\tau\in\mathcal{H}_F$.
\vskip.10in

Let $\xi$ be a periodic space-time white noise and $\rho:\mathbb{R}^4\rightarrow\mathbb{R}, \rho(t,x)=\rho_1(t)\rho_2(x)$ be a smooth compactly supported function integrating to $1$, set $\rho_\varepsilon(t,x)=\varepsilon^{-5}\rho(\frac{t}{\varepsilon^2},\frac{x}{\varepsilon})$.
Given the following  approximation $\xi_{\varepsilon,\vartheta}$ to  $\xi$,
 there is a canonical way of lifting it to an admissible  model $(\Pi^{(\varepsilon,\vartheta)}, \Gamma^{(\varepsilon,\vartheta)})$ as follows.
We set for $k\in\mathbb{Z}$,
$$\xi_\varepsilon=\rho_\varepsilon*\xi,\quad \xi_{\varepsilon,\vartheta}(t,x)=\frac{1}{\vartheta}\int_{k\vartheta}^{(k+1)\vartheta}\xi_\varepsilon(u,x)du, \quad t\in(k\vartheta,(k+1)\vartheta],$$

$$(\Pi_x^{(\varepsilon,\vartheta)}\Xi)(z)=\xi_{\varepsilon,\vartheta}(z), (\Pi_x^{({\varepsilon,\vartheta})}X^k)(z)=(z-x)^k$$
and
recursively define
$$(\Pi_x^{({\varepsilon,\vartheta})}\tau\bar{\tau})(z)=(\Pi_x^{({\varepsilon,\vartheta})}\tau)(z)(\Pi_x^{({\varepsilon,\vartheta})}\bar{\tau})(z),$$
and
\begin{equation}(\Pi_x^{({\varepsilon,\vartheta})}\mathcal{I}\tau)(z)=\int K(z-z_1)(\Pi_x^{({\varepsilon,\vartheta})}\tau)(z_1)dz_1+\sum_{l}\frac{(z-x)^l}{l!}f_x^{({\varepsilon,\vartheta})}(P_+\mathcal{I}_{l}\tau).\end{equation}
Here $f_x^{({\varepsilon,\vartheta})}(\mathcal{I}_{l}\tau)$ are defined by
\begin{equation}f_x^{({\varepsilon,\vartheta})}(\mathcal{I}_{l}\tau)=-\int D_1^lK(x-z_1)(\Pi_x^{({\varepsilon,\vartheta})}\tau)(z_1)dz_1.\end{equation}
Furthermore we impose $f_x^{({\varepsilon,\vartheta})}(X_i)=-x_i$, $f_x^{({\varepsilon,\vartheta})}(\tau\bar{\tau})=f_x^{(\varepsilon,\vartheta)}(\tau)f_x^{({\varepsilon,\vartheta})}(\bar{\tau})$ and extend this to all of $\mathcal{H}_{+}$ by linearity. Then define
\begin{equation}\Gamma_{xy}^{({\varepsilon,\vartheta})}=\Gamma_{f_x^{({\varepsilon,\vartheta})}}\circ (\Gamma_{f_y^{({\varepsilon,\vartheta})}})^{-1},\end{equation}
where $\Gamma_{f_x^{({\varepsilon,\vartheta})}}\tau:=(I\otimes {f_x^{({\varepsilon,\vartheta})}})\Delta \tau$ for $\tau\in\mathcal{H}_F$.

\vskip.10in
Then by [Hai14, Proposition 8.27] it is easy to check that $(\Pi^{({\varepsilon,\vartheta})},\Gamma^{({\varepsilon,\vartheta})})$ is an admissible  model for the regularity structure $\mathfrak{T}_F$ constructed in Section 2.1.

Now we give the following definition for the spaces of distributions $\mathcal{C}^\alpha_\mathfrak{s}$, $\alpha<0$, which is an extension of the definition of H\"{o}lder spaces to include $\alpha<0$.
\vskip.10in

\th{Definition 2.4} Let $\eta\in\mathcal{S}'$ and  $\alpha<0$. We say that $\eta\in \mathcal{C}^\alpha_\mathfrak{s}$ if the bound
$$|\eta(\varphi_z^\lambda)|\lesssim\lambda^\alpha,$$
holds uniformly over all $\lambda\in (0,1]$, all $\varphi\in\mathcal{B}_\alpha$ and locally uniformly over $z\in\mathbb{R}^{d+1}$.
\vskip.10in
For every compact set $\mathfrak{R}\subset \mathbb{R}^{d+1}$, we will denote by $\|\eta\|_{\alpha;\mathfrak{R}}$ the seminorm given by
$$\|\eta\|_{\alpha;\mathfrak{R}}:=\sup_{z\in\mathfrak{R}}\sup_{\varphi\in\mathcal{B}_\alpha}\sup_{\lambda\leq 1}\lambda^{-\alpha}|\eta(\varphi_z^\lambda)|.$$
We also write $\|\cdot\|_\alpha$ for the same expression with $\mathfrak{R}=\mathbb{R}^{d+1}$.
In the following we also use $\mathcal{C}^\alpha$ to denote $\mathcal{C}^\alpha_{\bar{\mathfrak{s}}}$ on $\mathbb{R}^d$ for the scaling $\bar{\mathfrak{s}}:=({\mathfrak{s}}_1,...,{\mathfrak{s}}_d)$.
\vskip.10in
We also have the following definition of spaces of  modelled distributions, which are the  H\"{o}lder spaces on the regularity structure.
Set $\mathfrak{P}=\{(t,x):t=0\}$.
Given a subset $\mathfrak{R}\subset \mathbb{R}^{d+1}$ we  denote by $\mathfrak{R}_\mathfrak{P}$ the set
$$\mathfrak{R}_\mathfrak{P}=\{(z,\bar{z})\in (\mathfrak{R}\setminus \mathfrak{P})^2:z\neq \bar{z}\textrm{ and } \|z-\bar{z}\|_\mathfrak{s}\leq |t|^{\frac{1}{\mathfrak{s}_0}}\wedge|\bar{t}|^{\frac{1}{\mathfrak{s}_0}}\wedge1\},$$
where $z=(t,x), \bar{z}=(\bar{t},\bar{x})$.

\vskip.10in
\th{Definition 2.5} Let $(\Pi,\Gamma)$ be a model for the regularity structure $\mathfrak{T}_F$  and $\mathfrak{P}$ as above. Then for any $\gamma>0$ and $\eta\in\mathbb{R}$, the space $\mathcal{D}^{\gamma,\eta}$  consists of all functions $f:\mathbb{R}^{d+1}\setminus \mathfrak{P}\rightarrow \bigoplus_{\alpha<\gamma}T_\alpha$ such that for every compact set $\mathfrak{R}\subset \mathbb{R}^{d+1}$ one has
$$|||f|||_{\gamma,\eta;\mathfrak{R}}:=\sup_{z\in \mathfrak{R}\setminus \mathfrak{P}}\sup_{l<\gamma}\frac{\|f(z)\|_l}{|t|^{\frac{\eta-l}{\mathfrak{s}_0}\wedge0}}+\sup_{(z,\bar{z})\in\mathfrak{R}_\mathfrak{P}}
\sup_{l<\gamma}\frac{\|f(z)-\Gamma_{z\bar{z}}f(\bar{z})\|_l}{\|z-\bar{z}\|^{\gamma-l}_\mathfrak{s}(|t|\wedge|\bar{t}|)^{\frac{\eta-\gamma}{\mathfrak{s}_0}}}<\infty.$$
Here  we wrote $\|\tau\|_l$ for the norm of the component of $\tau$ in $T_l$ and also used $t$ and $\bar{t}$ as shorthands for the time
components of the space-time points $z$ and $\bar{z}$.
\vskip.10in

For $f\in \mathcal{D}^{\gamma,\eta}$ and $\bar{f}\in \bar{\mathcal{D}}^{\gamma,\eta}$ (denoting by $\bar{\mathcal{D}}^{\gamma,\eta}$ the space built over another model $(\bar{\Pi},\bar{\Gamma})$), we also set
$$|||f;\bar{f}|||_{\gamma,\eta;\mathfrak{R}}:=\sup_{z\in \mathfrak{R}\setminus \mathfrak{P}}\sup_{l<\gamma}\frac{\|f(z)-\bar{f}(z)\|_l}{|t|^{\frac{\eta-l}{\mathfrak{s}_0}\wedge0}}+\sup_{(z,\bar{z})
\in\mathfrak{R}_\mathfrak{P}}\sup_{l<\gamma}\frac{\|f(z)-\bar{f}(z)-\Gamma_{z\bar{z}}f(\bar{z})+
\bar{\Gamma}_{z\bar{z}}\bar{f}(\bar{z})\|_l}{\|z-\bar{z}\|^{\gamma-l}_\mathfrak{s}(|t|\wedge|\bar{t}|)^{\frac{\eta-\gamma}{\mathfrak{s}_0}}},$$
which gives a natural
distance between elements $f\in \mathcal{D}^{\gamma,\eta}$ and $\bar{f}\in \bar{\mathcal{D}}^{\gamma,\eta}$.
\vskip.10in

Given a regularity structure, we say that a subspace $V\subset T$ is a sector of regularity $\alpha$ if it is invariant under the action of the structure group $G$ and it can be written as $V=\oplus_{\beta\in A}V_\beta$ with $V_\beta\subset T_\beta$, and $V_\beta=\{0\}$ for $\beta<\alpha$. We will use $\mathcal{D}^{\gamma,\eta}(V)$ to denote all functions in $\mathcal{D}^{\gamma,\eta}$ taking values in $V$.

Under suitable regularity assumptions, we can reconstruct from a given
modelled distribution $f$, a  distribution $\mathcal{R}f$ in the real world which "looks like $\Pi_xf(x)$ near $x$". This result, which defines the so-called reconstruction operator, is one of the most fundamental results in
 the theory of regularity structures.

\vskip.10in
\th{Theorem 2.6} (cf. [Hai14, Proposition 6.9]) Given a regularity structure and a model $(\Pi,\Gamma)$. Let $f\in \mathcal{D}^{\gamma,\eta}(V)$ for some sector $V$ of regularity $\alpha\leq0$, some $\gamma>0$, and some $\eta\leq \gamma$. Then provided that $\alpha\wedge\eta>-\mathfrak{s}_0$, there exists a unique distribution $\mathcal{R}f\in \mathcal{C}^{\eta\wedge\alpha}_\mathfrak{s}$ such that
$$|(\mathcal{R}f-\Pi_zf(z))(\varphi_z^\lambda)|\lesssim\lambda^\gamma,$$
holds uniformly over $\lambda\in (0,1]$ and $\varphi\in\mathcal{B}_r$ with $\varphi_z^\lambda$ compactly supported away from $\mathfrak{P}$ and locally uniformly over $z\in\mathbb{R}^{d+1}$.

Moreover, $(\Pi,\Gamma, f)\rightarrow\mathcal{R}f$ is jointly (locally) Lipschitz continuous with respect to the metric for $(\Pi,\Gamma)$ and $f$ defined in Definitions 2.2 and 2.5.
\vskip.10in

\subsection{Abstract fixed point problem}

We reformulate (1.1) as a fixed point problem in $\mathcal{D}^{\gamma,\eta}$ for suitable $\gamma$ and $\eta$. By Duhamel's formula, (1.1) is equivalent for smooth $\xi$ to the integral equation
$$u=G*((\xi-u^3)1_{t>0})+Gu_0.$$
Here, $G$ denotes the heat kernel, $*$ denotes space-time convolution, and $Gu_0$ denotes
the solution to the heat equation with initial condition $u_0$.
In order to interpret this equation as an identity in  $\mathcal{D}^{\gamma,\eta}$, we need  the following result from [Hai14, Proposition 6.16].
\vskip.10in
\th{Theorem 2.7} Let $\mathfrak{T}_F=(A,\mathcal{H}_F,G_F)$ be the regularity structure constructed  above and $(\Pi,\Gamma)$ be an admissible  model for $\mathfrak{T}_F$. Let $\gamma>0$, $\eta\leq \gamma$ and $\mathcal{I}$  act on some sector $V$ of regularity $\alpha\leq0$.  Then provided that $\alpha\wedge\eta>-2$, $\gamma+2, \eta+2$ not in $\mathbb{N}$, there exists a continuous linear operator $\mathcal{K}_\gamma:\mathcal{D}^{\gamma,\eta}(V)\rightarrow\mathcal{D}^{{\gamma'},{\eta'}}$ with ${\gamma'}=\gamma+2$ and ${\eta'}=(\eta\wedge\alpha)+2$,
such that
$$\mathcal{R}\mathcal{K}_\gamma f=K*\mathcal{R}f,$$
holds for $f\in \mathcal{D}^{\gamma,\eta}(V)$.
\vskip.10in
In the following we will only consider (1.1) with periodic boundary conditions. By the theory of regularity structures proposed in [Hai14] we can define translation maps and use it to define the periodic modelled distribution. Here the fundamental domain of the translation maps is compact.  We use the notations  $O_T=(-\infty,T]\times \mathbb{R}^d$ and use $|||\cdot|||_{\gamma,\eta;T}$ as a short hand for $|||\cdot|||_{\gamma,\eta;O_T}$.
Moreover, we have that for   $\gamma,\eta,\gamma',\eta'$ in Theorem 2.7 and some $\theta>0$
$$|||\mathcal{K}_\gamma 1_{t>0}f|||_{\gamma',\eta';T}\lesssim T^\theta|||f|||_{\gamma,\eta;T}.$$

Now we reformulate the fixed point map as
\begin{equation}\aligned v=&(\mathcal{K}_{\bar{\gamma}}+R_\gamma\mathcal{R})(\mathbf{1}_{t>0}\Xi),
\\u=&-(\mathcal{K}_{\bar{\gamma}}+R_\gamma\mathcal{R})(\mathbf{1}_{t>0} u^{3})+v+\mathcal{G}u_0
.\endaligned\end{equation}
Here $\mathbf{1}_{t>0}(t,x)=1$ for $t>0$ and $\mathbf{1}_{t>0}(t,x)=0$ otherwise,  and for the smooth function $R=G-K$,
 $$R_\gamma:\mathcal{C}^\alpha_\mathfrak{s}\rightarrow\mathcal{D}^{\gamma,\eta}, (R_\gamma f)(z)=\sum_{|k|_\mathfrak{s}<\gamma}\frac{X^k}{k!}\int D^kR(z-\bar{z})f(\bar{z})d\bar{z},$$
 $$\mathcal{G}u_0=\sum_{|k|_\mathfrak{s}<\gamma}\frac{X^k}{k!} D^k(Gu_0)(z),$$
where $\gamma,\bar{\gamma}$ will be chosen below and we define $\mathcal{R}(\mathbf{1}_{t>0}\Xi)$ as the distribution $\xi\textbf{1}_{t>0}$.

We consider the second equation in (2.7): Define
$$V:=\mathcal{I}(\mathcal{F}_F) \oplus\bar{T}.$$

\vskip.10in
Now for   $u_0\in \mathcal{C}^\eta(\mathbb{R}^3)$, $\eta$ not in $\mathbb{N}$, periodic,  [Hai14, Lemma 7.5] implies that $\mathcal{G}u_0\in \mathcal{D}^{\gamma,\eta}$ for $\gamma>(\eta\vee0)$.

We define for any $\beta<0$ and any compact set $\mathfrak{R}$ the norm
$$|\xi|_{\beta;\mathfrak{R}}=\sup_{s\in\mathbb{R}}\|\xi1_{t\geq s}\|_{\beta;\mathfrak{R}},$$
and we denote by $\bar{\mathcal{C}}_\mathfrak{s}^\beta$ the intersections of the completions of smooth functions under $|\cdot|_{\beta;\mathfrak{R}}$ for all compact sets $\mathfrak{R}$.  By [Hai14, Proposition 9.5] we know that for every $\alpha\in (-3,-\frac{5}{2})$, the space-time white noise $\xi$ belongs to $\bar{\mathcal{C}}_\mathfrak{s}^\alpha$ almost surely and  $K*\xi\in C(\mathbb{R},\mathcal{C}^{\alpha+2}(\mathbb{R}^3))$ almost surely. With these notations at hand, we recall the following results from [Hai14].
\vskip.10in
\th{Proposition 2.8}( [Hai14, Proposition 9.8]) Let $\mathfrak{T}_F$ be the regularity structure  associated to ($\Phi^4$) with $\alpha\in (-\frac{18}{7},-\frac{5}{2})$. Let $\eta\in (-\frac{2}{3},\alpha+2)$, $\gamma>|2\alpha+4|$, $\bar{\gamma}=\gamma+2\alpha+4$ and let $Z=(\Pi,\Gamma)$ be an admissible model for $\mathfrak{T}_F$  with the additional properties that  $\xi:=\mathcal{R}\Xi$ belongs to $\bar{\mathcal{C}}_\mathfrak{s}^\alpha$ and that $K*\xi\in C(\mathbb{R},\mathcal{C}^{\eta})$. Then there exists a maximal solution $\mathcal{S}^L(u_0,Z)\in\mathcal{D}^{\gamma,\eta}(V)$ to the equation (2.7).

Furthermore, let  $ T^L(u_0, Z)\in \mathbb{R}^+\cup\{+\infty\}$ be the first time such
that  $\|(\mathcal{R}\mathcal{S}^L(u_0,Z))(t,\cdot)\|_\eta\geq L$ and set $O=[-1,2]\times \mathbb{R}^d$. Then, for every $\varepsilon> 0$ and $C > 0$
there exists $\delta> 0$ such that, setting $T =1\wedge T^L(u_0, Z)\wedge T^L(\bar{u}_0,\bar{Z})$, one has the bound
$|||\mathcal{S}^L(u_0, Z)-\mathcal{S}^L(\bar{u}_0,\bar{Z})|||_{\gamma,\eta;T} \leq \varepsilon$,
for all $u_0, \bar{u}_0, Z,\bar{Z}$ satisfying $|||Z|||_{\gamma;O}\leq C, |||\bar{Z}|||_{\gamma;O}\leq C, \|u_0\|_\eta\leq L/2, \|\bar{u}_0\|_\eta\leq L/2,
\|u_0-\bar{u}_0\|_\eta\leq \delta$, and $|||Z;\bar{Z}|||_{\gamma;O}\leq \delta$,
and satisfying the bounds
$|\xi|_{\alpha;O} + |\bar{\xi}|_{\alpha;O}\leq C , \sup_{t\in[0,1]}
\|(K*\xi)(t,\cdot)\|_\eta+\sup_{t\in[0,1]}\|(K*\bar{\xi})(t,\cdot)\|_\eta\leq C,$ as well as
$$|\xi-\bar{\xi}|_{\alpha;O}\leq \delta , \sup_{t\in[0,1]}
\|(K*\xi)(t,\cdot)-(K*\bar{\xi})(t,\cdot)\|_\eta\leq \delta,$$
Here we have set $\bar{\xi}=\bar{\mathcal{R}}\Xi$, where $\bar{\mathcal{R}}$ is the reconstruction operator associated to $\bar{Z}$.
\vskip.10in

\section{Renormalisation procedure and main result}
In Section 2 we have constructed a  model associated with $\xi_{\varepsilon,\vartheta}$ and in this section we will prove the  convergence result required in Proposition 2.8, which at last implies Theorem 1.1. As we mentioned in the introduction,  the sequence of models does not converge to  a limit. We have to renormalise the model into some converging renormalised model.

\subsection{Renormalised model}

In this subsection we  renormalise the model and prove that it is also an admissible model for the regularity structure $\mathfrak{T}_F$ associated with the dynamical $\Phi^4_3$ model. In our case we should subtract some functions (denoted by $C_1^{(\varepsilon,\vartheta)}, C_2^{(\varepsilon,\vartheta)}$) depending on $t$ in the renormalisations (see the proof of Theorem 3.7), which cannot be written as the sum of diverging constants and converging functions as explained in Remark 1.2.  This is the main difference from the case in [Hai14], where all the terms being subtracted in the renormalisations are constants.  To prove that the  renormalised model is admissible in our case, we define a bigger regularity structure $\mathfrak{T}^1$  including the original regularity structure $\mathfrak{T}_F$ and two symbols $\mathbf{C}_1, \mathbf{C}_2$, where $\mathbf{C}_1$ and $\mathbf{C}_2$ represent $C_1^{(\varepsilon,\vartheta)}, C_2^{(\varepsilon,\vartheta)}$ in the regularity structure, respectively. We build a model  for $\mathfrak{T}^1$ and use it to prove the desired result.

First, we construct the regularity structure $\mathfrak{T}^1$.  Define a set $\mathcal{F}^1$ by postulating that $\{\mathbf{1},\Xi,X_j,$ $\mathbf{C}_1, \mathbf{C}_2\}\subset \mathcal{F}^1$ and whenever $\tau,\bar{\tau}\in\mathcal{F}^1$, we have
$\tau\bar{\tau}\in\mathcal{F}^1$ and $\mathcal{I}_k(\tau)\in\mathcal{F}^1$; define $\mathcal{F}_+^1$ as the set of all elements $\tau\in\mathcal{F}^1$ such that either $\tau=\mathbf{1}$ or $|\tau|_\mathfrak{s}>0$ and such that, whenever $\tau$ can be written as $\tau=\tau_1\tau_2$ we have either $\tau_i=\mathbf{1}$ or $|\tau_i|_\mathfrak{s}>0$; $\mathcal{H}^1, \mathcal{H}_+^1$ denote the sets of finite linear combinations of all elements in $\mathcal{F}^1, \mathcal{F}_+^1$, respectively.
Here for each $\tau\in \mathcal{F}^1$ a weight $|\tau|_{\mathfrak{s}}$  is defined as in Section 2 and by setting
$|\mathbf{C}_1|_\mathfrak{s}=|\mathbf{C}_2|_\mathfrak{s}=-\delta_0$ with $4\alpha+10<-\delta_0<0$.
The reason for $4\alpha+10<-\delta_0<0$ is  to make sure that the homogeneity of $\mathbf{C}_2$ is bigger than $\mathcal{I}(\Psi^2)\Psi^2$.

Recall that $\mathfrak{M}_F=\{\Xi, U^n:n\leq 3\}.$
We define the sets $\mathcal{W}_n^1, \mathcal{U}_n^1$ for $n\geq0$ recursively by
$$\mathcal{W}_0^1=\mathcal{U}_0^1=\emptyset,$$
$$\mathcal{W}_n^1=\mathcal{W}_{n-1}^1\cup\bigcup_{\mathcal{Q}\in\mathfrak{M}_F}\mathcal{Q}(\mathcal{U}_{n-1}^1,\Xi),$$
$$\mathcal{U}_n^1=\{X^k, \mathbf{C}_1, \mathbf{C}_2\}\cup\{\mathcal{I}(\tau):\tau\in\mathcal{W}_n\}$$
and
$$\mathcal{F}_F^1:=\bigcup_{n\geq0}(\mathcal{W}_n^1\cup\mathcal{U}_n^1).$$

We denote by $\mathcal{H}_F^1$ the set of all finite linear combinations of elements in $\mathcal{F}_F^1$ and denote by $\mathcal{F}_F^{1,+}$ the set of those basis vectors $\bar{\tau}\in\mathcal{F}_+^1$ that can be written as $\bar{\tau}=X^{l_0}\Pi_i\mathcal{I}_{l_i}\tau_i$ for some multiindices $l_i$ and some elements $\tau_i\in\mathcal{F}_F^1$. Denote $\mathcal{H}_F^{1,+}$ the set of all finite linear combinations of elements in $\mathcal{F}_F^{1,+}$.
\vskip.10in

Now we  construct the structure group $G_F^1$. We  define the operators from $\mathcal{H}^1$ to $\mathcal{H}^{1}\otimes \mathcal{H}^{1}_+$ and from $\mathcal{H}^1_+$ to $\mathcal{H}^{1}_+\otimes \mathcal{H}^{1}_+$ as $\Delta, \Delta^+$ in Section 2. We still use $\Delta, \Delta^+$ to denote them for  notational simplicity. $\Delta$ on $\mathbf{1}$, $X_i$, $\Xi$ and $\Delta^+$ on $\mathbf{1}$, $X_i$  can be defined as in Section 2.  Define
$$\Delta \mathbf{C}_1=\mathbf{C}_1\otimes \mathbf{1},\quad \Delta \mathbf{C}_2=\mathbf{C}_2\otimes \mathbf{1}.$$
For all other terms $\Delta, \Delta^+$  can also be defined recursively as in Section 2.

 By using the theory of regularity structures (see [Hai14, Section 8]) we can  define a structure group $G_F^1$ of linear operators acting on  $\mathcal{H}_F^1$ satisfying Definition 2.1 as follows: For  $g\in \mathcal{H}_+^{1,*}$, the dual of $\mathcal{H}_+^1$, satisfying $g(\tau \bar{\tau})=g(\tau)g(\bar{\tau})$ for $\tau,\bar{\tau}\in\mathcal{H}_+^1$, define  $\Gamma_g:\mathcal{H}^1\rightarrow\mathcal{H}^1, \Gamma_g\tau=(I\otimes g)\Delta\tau$.  By [Hai14, Theorem 8.24] we construct the following regularity structure.

\vskip.10in
\th{Theorem 3.1} Let $T=\mathcal{H}_F^1$ with $T_\gamma=\langle\{\tau\in\mathcal{F}_F^1:|\tau|_\mathfrak{s}=\gamma\}\rangle$, $A^1=\{|\tau|_\mathfrak{s}:\tau\in\mathcal{F}_F^1\}$. Then $\mathfrak{T}^1=(A^1,\mathcal{H}_F^1,G_F^1)$ defines a regularity structure $\mathfrak{T}^1$.
\vskip.10in
We emphasize that we do not change the regularity structure associated with $\Phi^4$ in our case.  The introduction of $\mathfrak{T}^1$ is to prove that the renormalised model is an admissible model for $\mathfrak{T}_F$. In the following we extend the model $(\Pi^{({\varepsilon,\vartheta})},\Gamma^{({\varepsilon,\vartheta})})$ constructed in Section 2 to a model for $\mathfrak{T}^1$, which is used to construct the renormalised model. We still denote it by $(\Pi^{({\varepsilon,\vartheta})},\Gamma^{({\varepsilon,\vartheta})})$ for simplicity.

Given continuous functions $C_1^{(\varepsilon,\vartheta)}(t), C_2^{(\varepsilon,\vartheta)}(t)$, for $z=(t,y)$ we extend the models as follows:
$$(\Pi_x^{(\varepsilon,\vartheta)}\mathbf{C}_1)(z)=C_1^{(\varepsilon,\vartheta)}(t),\quad (\Pi_x^{(\varepsilon,\vartheta)}\mathbf{C}_2)(z)=C_2^{(\varepsilon,\vartheta)}(t),$$
and
recursively as in Section 2.2. Moreover, we extend $f_x^{(\varepsilon,\vartheta)}$ to all of $\mathcal{H}^{1,+}_F$ by linearity. Moreover,
$\Gamma_{xy}^{({\varepsilon,\vartheta})}$ is
still given by (2.6).

\vskip.10in
\th{Proposition 3.2}  $(\Pi^{({\varepsilon,\vartheta})},\Gamma^{({\varepsilon,\vartheta})})$ is a model for the regularity structure $\mathfrak{T}^1$ constructed in Theorem 3.1.

\proof Since $C_1^{(\varepsilon,\vartheta)}$, $C_2^{(\varepsilon,\vartheta)}$ are continuous functions,   a similar argument as in the proof of [Hai14, Proposition 8.27] implies the result.$\hfill\Box$

\vskip.10in

Now we introduce the following sets as in [Hai14, Section 9]. Define
$$\aligned\mathcal{F}_0:=\{&\mathbf{1},\Xi,\Psi, \Psi^2, \Psi^3, \Psi^2X_i, \mathcal{I}(\Psi^3)\Psi, \mathcal{I}(\Psi^3)\Psi^2,
\\&\mathcal{I}(\Psi^2)\Psi^2, \mathcal{I}(\Psi^2), \mathcal{I}(\Psi)\Psi, \mathcal{I}(\Psi)\Psi^2, X_i \},\endaligned$$
$$\mathcal{F}_*:=\{\Psi,\Psi^2,\Psi^3\},$$
where $\Psi=\mathcal{I}(\Xi)$ and the index $i$ corresponds to any of the three spatial directions.

Then $\mathcal{F}_0\subset \mathcal{F}_F$ contains every $\tau\in\mathcal{F}_F$ with $|\tau|_\mathfrak{s}\leq0$ and for every $\tau\in\mathcal{F}_0$, $\Delta\tau\in\mathcal{H}_0\otimes\mathcal{H}_0^+$. Here $\mathcal{H}_0$ denotes the linear span of $\mathcal{F}_0$ and $\mathcal{H}_0^+$ denotes the linear span of the  elements in $\mathcal{F}_+$ of the form $X^k\prod_{i}\mathcal{I}_{l_i}\tau_i$ for some multiindices $k$ and $l_i$ such that $|\mathcal{I}_{l_i}\tau_i|_\mathfrak{s}>0$ and $\tau_i\in\mathcal{F}_*$.

With these notations at hand, we construct a linear map $M$ from $\mathcal{H}_0$ to $\mathcal{H}^1_F$ by
\begin{equation}\aligned M\Psi^2=&\Psi^2-\mathbf{C}_1,\\M(\Psi^2X_i)=&\Psi^2X_i-\mathbf{C}_1X_i,\\M\Psi^3=&\Psi^3-3\mathbf{C}_1\Psi,\\M\mathcal{I}(\Psi^2)=&\mathcal{I}(\Psi^2)
-\mathcal{I}(\mathbf{C}_1),\\M(\mathcal{I}(\Psi^2)\Psi^2)=&(\mathcal{I}(\Psi^2)
-\mathcal{I}(\mathbf{C}_1))
(\Psi^2-\mathbf{C}_1)-\mathbf{C}_2,\\M(\mathcal{I}(\Psi^3)\Psi)=&(\mathcal{I}(\Psi^3)-3\mathcal{I}(\mathbf{C}_1\Psi))\Psi,\\M(\mathcal{I}(\Psi^3)\Psi^2)=&(\mathcal{I}(\Psi^3)-3\mathcal{I}(\mathbf{C}_1\Psi))(\Psi^2-\mathbf{C}_1)
-3\mathbf{C}_2\Psi,\\M(\mathcal{I}(\Psi)\Psi^2)=&\mathcal{I}(\Psi)(\Psi^2
-\mathbf{C}_1),\endaligned\end{equation}
as well as $M\tau=\tau$ for the remaining basis elements $\tau\in\mathcal{F}_0$. In our case, $\mathbf{C}_i, i=1,2,$ are not in $\bar{T}$ and hence, $\mathcal{I}(\mathbf{C}_i)\neq0$.  Now similarly as in [Hai14, Section 8] we introduce the following linear maps which are used to construct the renormalised model.  Define
a linear map $\Delta^M:\mathcal{H}_0\rightarrow\mathcal{H}_F^1\times \mathcal{H}_F^{1,+}$ by
$$\Delta^M\tau=(M\tau)\otimes \mathbf{1},$$
for those elements $\tau\in\mathcal{F}_0$  not containing a factor $\mathcal{I}(\Psi^2)$ or $\mathcal{I}(\Psi^3)$. For the remaining  elements, we define
$$\Delta^M\mathcal{I}(\Psi^2)=(M(\mathcal{I}(\Psi^2)))\otimes \mathbf{1}+ X_i\otimes \mathcal{I}_i(\mathbf{C}_1),$$
$$\Delta^M\mathcal{I}(\Psi^2)\Psi^2=(M(\mathcal{I}(\Psi^2)\Psi^2))\otimes \mathbf{1}+ (\Psi^2-\mathbf{C}_1)X_i\otimes \mathcal{I}_i(\mathbf{C}_1),$$
$$\Delta^M\mathcal{I}(\Psi^3)\Psi=(M(\mathcal{I}(\Psi^3)\Psi))\otimes \mathbf{1}+3\Psi X_i\otimes \mathcal{I}_i(\mathbf{C}_1\Psi),$$
$$\Delta^M\mathcal{I}(\Psi^3)\Psi^2=(M(\mathcal{I}(\Psi^3)\Psi^2))\otimes \mathbf{1}+3(\Psi^2-\mathbf{C}_1) X_i\otimes \mathcal{I}_i(\mathbf{C}_1\Psi).$$
Moreover, we introduce a  linear map  $\hat{M}:\mathcal{H}_0^+\rightarrow\mathcal{H}_F^{1,+}$, which is a multiplicative morphism and  leaves $X^k$ invariant, and
$$\hat{M}\mathcal{I}(\Psi^n)=\mathcal{I}(M\Psi^n),\quad \hat{M}\mathcal{I}_i(\Psi)=\mathcal{I}_i(\Psi).$$
Then we can easily check that
\begin{equation}\hat{M}\mathcal{I}_k=\mathcal{M}(\mathcal{I}_k\otimes I)\Delta^M,\end{equation}
\begin{equation}(I\otimes \mathcal{M})(\Delta\otimes I)\Delta^M=(M\otimes \hat{M})\Delta.\end{equation}
Here $\mathcal{M}:\mathcal{H}_F^{1,+}\times \mathcal{H}_F^{1,+}\rightarrow\mathcal{H}_F^{1,+}$ denotes the multiplication map.

Furthermore, define a linear multiplicative morphism:
 $\hat{\Delta}^M:\mathcal{H}_0^+\rightarrow\mathcal{H}_F^{1,+}\times\mathcal{H}_F^{1,+}$ by
 $$\hat{\Delta}^M X^k=X^k\otimes \mathbf{1},$$
and $$\aligned\hat{\Delta}^M\mathcal{I}(\Psi^n)=&\mathcal{I}(M\Psi^n)\otimes \mathbf{1}+3\delta_{n3}(X_i\otimes \mathcal{I}_i(\mathbf{C}_1\Psi)-X_i\mathcal{I}_i(\mathbf{C}_1\Psi)\otimes \mathbf{1})\\&+\delta_{n2}(X_i\otimes\mathcal{I}_i(\mathbf{C}_1)-X_i\mathcal{I}_i(\mathbf{C}_1)\otimes \mathbf{1}).\endaligned$$
Then we can easily check that
$$(\mathcal{A}\hat{M}\mathcal{A}\otimes \hat{M})\Delta^+=(I\otimes \mathcal{M})(\Delta^+\otimes I)\hat{\Delta}^M.$$
Here $\mathcal{A}$ is as given in [Hai14, Section 8] for the regularity structure $\mathfrak{T}^1$.

Now we give the renormalised model by using the above maps: Define for $\tau\in\mathcal{H}_0, \tau_1\in\mathcal{H}_0^+$,
 $$\Pi_x^{M,\varepsilon,\vartheta}\tau=(\Pi_x^{({\varepsilon,\vartheta})}\otimes f_x^{({\varepsilon,\vartheta})})\Delta^M\tau,\quad f_x^{M,\varepsilon,\vartheta}\tau_1=f_x^{({\varepsilon,\vartheta})}\hat{M}\tau_1,$$
 and define $\Gamma^M_{xy}=(F_x^M)^{-1}\circ F_y^M$ with $F_x^M:=(I\otimes f_x^M)\Delta$. Then by a similar argument as in the proof of [Hai14, Theorem 8.44] we have the following result.

\vskip.10in
  \th{Proposition 3.3} $(\Pi^M,\Gamma^M)$ is an admissible model for $\mathfrak{T}_F$ on $\mathcal{H}_0$. Furthermore, it extends uniquely to an admissible model for all of $\mathfrak{T}_F$.

 \proof By the definition of $\Pi^M$ and the expression for $\Delta^M$ we know that $(\Pi_x^M\tau)(\varphi_x^\lambda)$ can be written as a finite linear combination of terms of the type $(\Pi_x\bar{\tau})(\varphi_x^\lambda)$ with $|\bar{\tau}|_\mathfrak{s}\geq|\tau|_\mathfrak{s}$ and $\bar{\tau}\in\mathcal{H}^1_F$. Then  Proposition 3.2 implies the required scaling as a function of $\lambda$.

Define $\gamma_{xy}:=(f_x\mathcal{A}\otimes f_y)\Delta^+$ and we have
$\Gamma_{xy}=(I\otimes\gamma_{xy})\Delta$.
 Since $(\Pi^{({\varepsilon,\vartheta})},\Gamma^{({\varepsilon,\vartheta})})$ is a model for $\mathfrak{T}^1$, this implies that for $\tau\in\mathcal{H}_F^{1,+}$
 $$|\gamma_{xy}\tau|\lesssim\|x-y\|_\mathfrak{s}^{|\tau|_\mathfrak{s}}.$$
 Since $\Gamma^M_{xy}=(I\otimes\gamma_{xy}^M)\Delta$ with $\gamma_{xy}^M=(\gamma_{xy}\otimes f_y)\hat{\Delta}^M$ and $\hat{\Delta}^M\tau=\tau\otimes \mathbf{1}+\sum \tau^1\otimes \tau^2$ with $|\tau^1|_{\mathfrak{s}}>|\tau|_{\mathfrak{s}}$,  it follows from the expression of $\hat{\Delta}^M$ that for $\tau\in\mathcal{H}_0^+$
 $$|\gamma_{xy}^M\tau|\lesssim\|x-y\|_\mathfrak{s}^{|\tau|_\mathfrak{s}}.$$
 Thus $(\Pi^M,\Gamma^M)$ is  a model on $\mathcal{H}_0$. By (3.2), (3.3) and similar arguments as in [Hai14, Section 8] we know that $(\Pi^M,\Gamma^M)$ is also an admissible model on  $\mathcal{H}_0$.  Finally  applying [Hai14, Theorem 5.14,  Proposition 3.31]
$(\Pi^M,\Gamma^M)$ can be extended uniquely to all of $\mathfrak{T}_F$.$\hfill\Box$
\vskip.10in

\th{Remark 3.4} (i) It is a little different from the case in [Hai14] to construct the renormalised model. In [Hai14] the renormalised map $M$ is a linear map from $\mathcal{H}_0$ to $\mathcal{H}_0$, which is enough for the construction of the renormalised model. In our case we have to subtract some functions $C_1, C_2$ to make the diverging terms converge in some sense. As we explained at the beginning of the section, we construct a new regularity structure $\mathfrak{T}^1$  including $\mathbf{C}_1, \mathbf{C}_2$ which represent the functions $C_1$ and $C_2$, respectively. The renormalised map $M$ is a linear map from $\mathcal{H}_0$ to  $\mathcal{H}^1_F$, which does not belong to the renormalisation group defined in [Hai14, Definition 8.41]. However, we could still use it to define $\hat{M}, \Delta^M, \hat{\Delta}^M$ and construct  the renormalised model on $\mathcal{H}_F$. We emphasize that the renormalised model $(\Pi^M,\Gamma^M)$ is associated with the regularity structure $\mathfrak{T}_F$. Below we still consider the  regularity structure $\mathfrak{T}_F$. $\mathfrak{T}^1$ is a tool to prove that the renormalised model is an admissible model for $\mathfrak{T}_F$.

(ii)  In fact, we can also define the renormalised model for the bigger regularity structure $\mathfrak{T}^1$ and apply directly the results in [Hai14, Section 8] to conclude that  the renormalised model is an admissible model for $\mathfrak{T}^1$, which is also the required renormalised model when restricted on $\mathfrak{T}_F$. For this argument we need to define the corresponding $\mathcal{F}_0$ for $\mathfrak{T}^1$, which is a little bit complicated. Therefore, we rather use the above proof, because it appears to be simpler.
\vskip.10in

\subsection{Renormalised solutions}

Let $u_{\varepsilon,\vartheta}=\bar{\mathcal{S}}^L(u_0,\xi_{\varepsilon,\vartheta})$ denote the classical solution map to the equation
$$\partial_tu_{\varepsilon,\vartheta}=\Delta u_{\varepsilon,\vartheta}-u_{\varepsilon,\vartheta}^3+\xi_{\varepsilon,\vartheta},\quad u_{\varepsilon,\vartheta}(0)=u_0.$$
 Here $u_0\in \mathcal{C}^\eta(\mathbb{T}^3)$. The renormalised map $\bar{\mathcal{S}}^L_M(u_0,\xi_{\varepsilon,\vartheta})$ is given by the classical solution map to the equation
$$\partial_tu_{\varepsilon,\vartheta}=\Delta u_{\varepsilon,\vartheta}+(3C_1^{(\varepsilon,\vartheta)}-9C_2^{(\varepsilon,\vartheta)})u_{\varepsilon,\vartheta}-u_{\varepsilon,\vartheta}^3
+\xi_{\varepsilon,\vartheta},\quad u_{\varepsilon,\vartheta}(0)=u_0.$$
By the same argument as in the proof of [Hai14, Proposition 9.10] we obtain the following result:
\vskip.10in

\th{Proposition 3.5} Let  $Z_{\varepsilon,\vartheta}=(\Pi^{(\varepsilon,\vartheta)},\Gamma^{(\varepsilon,\vartheta)})$ denote the model given in Section 2, and
 $Z_{\varepsilon,\vartheta}^M=(\Pi^{M,\varepsilon,\vartheta},\Gamma^{M,\varepsilon,\vartheta})$ the renormalised model in Propostion 3.3. Then for every  $u_0\in\mathcal{C}^\eta(\mathbb{T}^3)$ one has the identities
$$\mathcal{R}\mathcal{S}^L(u_0,Z_{\varepsilon,\vartheta})=\bar{\mathcal{S}}^L(u_0,\xi_{\varepsilon,\vartheta}),\quad \mathcal{R}\mathcal{S}^L(u_0,Z_{\varepsilon,\vartheta}^M)=\bar{\mathcal{S}}^L_M(u_0,\xi_{\varepsilon,\vartheta}).$$
Here $\mathcal{S}^L(u_0,Z_{\varepsilon,\vartheta})$ and $\mathcal{S}^L(u_0,Z_{\varepsilon,\vartheta}^M)$ are the solutions obtained in Proposition 2.8.
\vskip.10in

\subsection{Proof of the main result}

In this subsection we prove Theorem 1.1. We first prove the required convergence in Proposition 2.8 for $\xi_{\varepsilon,\vartheta}$ and $K*\xi_{\varepsilon,\vartheta}$. Our argument essentially follows [Hai14, Proposition 9.5].

\vskip.10in
\th{Proposition 3.6} Let $\xi$ be white noise on $\mathbb{R}\times \mathbb{T}^3$, which we extend periodically to $\mathbb{R}^4$,  and define $\xi_{\varepsilon,\vartheta}$ as in Subsection 2.2. Then for every compact set $\mathfrak{R}\subset \mathbb{R}^4$ and every $0<\kappa<-\alpha-\frac{5}{2}$
we have
\begin{equation}\textbf{E}|\xi_{\varepsilon,\vartheta} -\xi|_{\alpha;\mathfrak{R}}\lesssim \varepsilon^\kappa+\vartheta^\frac{\kappa}{2}.\end{equation}
Finally for every $0<\kappa<-\frac{2\alpha+5}{4}$, the bound
$$\textbf{E}\sup_{t\in[0,1]}\|K*\xi_{\varepsilon,\vartheta}(t,\cdot)-K*\xi(t,\cdot)\|_{\alpha+2}\lesssim\varepsilon^{2\kappa}+\vartheta^\kappa,$$
holds uniformly over $\varepsilon,\vartheta\in(0,1]$.

\proof
For any scaling $\mathfrak{s}$ of $\mathbb{R}^{4}$ and any $n\in\mathbb{Z}$, define
$$\Lambda^n_\mathfrak{s}=\{\sum_{j=0}^3 2^{-n\mathfrak{s}_j}k_je_j:k_j\in\mathbb{Z}\},$$
where  $e_j$ denotes the jth element of the canonical basis of $\mathbb{R}^4$.
We choose a wavelet basis $\{\psi_x^{n,\mathfrak{s}}=2^{-\frac{5}{2}n}\psi^{2^{-n}}_x,\varphi(\cdot-y), n\geq0, x\in \Lambda^n_\mathfrak{s}\cap\mathfrak{R}, y\in \Lambda^0_\mathfrak{s}\cap\mathfrak{R}, \psi\in \Psi\}$ as in [Hai14, Section 3.2] on $\mathbb{R}^4$. Writing $\Psi_\star=\Psi\cup\{\varphi\}$,
 we note that for every $p>1$, we have the bound
$$\aligned &\textbf{E}\|(\xi_{\varepsilon,\vartheta}-\xi)1_{t\in[0,s]}\|^{2p}_{\alpha;\mathfrak{R}}\\\leq& \sum_{\psi\in\Psi_\star}\sum_{n\geq0}\sum_{x\in\Lambda^n_\mathfrak{s}\cap\bar{\mathfrak{R}}}\textbf{E}2^{2\alpha np+|\mathfrak{s}|np}
|\langle (\xi_{\varepsilon,\vartheta}-\xi)1_{t\in[0,s]},\psi_{x}^{n,\mathfrak{s}}\rangle|^{2p}\\\lesssim&\sum_{\psi\in\Psi_\star}
\sum_{n\geq0}\sum_{x\in\Lambda^n_\mathfrak{s}\cap\bar{\mathfrak{R}}}2^{2\alpha np+|\mathfrak{s}|np}
(\textbf{E}|\langle (\xi_{\varepsilon,\vartheta}-\xi)1_{t\in[0,s]},\psi_{x}^{n,\mathfrak{s}}\rangle|^{2})^p.\endaligned$$
Here we wrote $\bar{\mathfrak{R}}$ for the $1$-fattening of $\mathfrak{R}$.
Since it has been obtained in [Hai14, Proposition 9.5] that
$$\aligned \textbf{E}|\langle (\xi_{\varepsilon}-\xi)1_{t\in[0,s]},\psi_{x}^{n,\mathfrak{s}}\rangle|^{2}\lesssim&1\wedge (2^{2n}s)\wedge(2^{2n}\varepsilon^2),\endaligned$$
 for $\xi_\varepsilon=\xi*\rho_\varepsilon$, it suffices to estimate $\textbf{E}|\langle (\xi_{\varepsilon,\vartheta}-\xi_\varepsilon)1_{t\in[0,s]},\psi_{x}^{n,\mathfrak{s}}\rangle|^{2}$.
By the definition of $\xi_{\varepsilon,\vartheta}$ we know that  $$\aligned &1_{t\in[0,s]}[\xi_{\varepsilon,\vartheta}(t)-\xi_\varepsilon(t)]\\=&\sum_{k=0}^{[\frac{s}{\vartheta}]-1}1_{t\in[k\vartheta,(k+1)\vartheta)}
\frac{1}{\vartheta}\int_{k\vartheta}^{(k+1)\vartheta}(\xi_{\varepsilon}(u)-\xi_\varepsilon (t))du+1_{t\in[[\frac{s}{\vartheta}]\vartheta,s]}\frac{1}{\vartheta}
\int_{[\frac{s}{\vartheta}]\vartheta}^{([\frac{s}{\vartheta}]+1)\vartheta}\xi_{\varepsilon}(u)du-1_{t\in[[\frac{s}{\vartheta}]\vartheta,s]}\xi_\varepsilon(t).
\\:=&I_1+I_2-I_3.\endaligned$$
In the following we  estimate $\textbf{E}|\langle I_i,\psi_{x}^{n,\mathfrak{s}}\rangle|^{2}$ for $i=1,2,3$ separately.
Since $\|\rho_\varepsilon*f\|_{L^2}\leq \|f\|_{L^2}$,  a straightforward calculation  yields that
$$\aligned \textbf{E}|\langle I_1,\psi_{x}^{n,\mathfrak{s}}\rangle|^{2}=&\textbf{E}|\int \sum_{k=0}^{[\frac{s}{\vartheta}]-1}\frac{1}{\vartheta}\int_{k\vartheta}^{(k+1)\vartheta}\int_{k\vartheta}^{(k+1)\vartheta}
(\psi_{x}^{n,\mathfrak{s}}(u,y)-\psi_{x}^{n,\mathfrak{s}}(t,y))\xi_\varepsilon(t,y)du dtdy|^2\\\lesssim&\int\int\bigg(\sum_{k=0}^{[\frac{s}{\vartheta}]-1}
\frac{1}{\vartheta}1_{t\in[k\vartheta,(k+1)\vartheta)}\int_{k\vartheta}^{(k+1)\vartheta}
[\psi_{x}^{n,\mathfrak{s}}(u,y)-\psi_{x}^{n,\mathfrak{s}}(t,y)]du\bigg)^2dtdy\\\lesssim&\frac{1}{\vartheta^2}\int\sum_{k=0}^{[\frac{s}{\vartheta}]-1}
\int_{k\vartheta}^{(k+1)\vartheta}\bigg(\int_{k\vartheta}^{(k+1)\vartheta}
\int_{k\vartheta}^{(k+1)\vartheta}|D_{\tilde{u}}\psi_{x}^{n,\mathfrak{s}}
(\tilde{u},y)|d\tilde{u}du\bigg)^2dtdy\\\lesssim&\vartheta^2\int_0^s\int |D_{\tilde{u}}\psi_{x}^{n,\mathfrak{s}}
(\tilde{u},y)|^2d\tilde{u}dy\\\lesssim&1\wedge (2^{2n}s)\wedge(2^{2n}\vartheta).\endaligned$$
Similarly,
$$\aligned\textbf{E}|\langle I_2,\psi_{x}^{n,\mathfrak{s}}\rangle|^{2}+\textbf{E}|\langle I_3,\psi_{x}^{n,\mathfrak{s}}\rangle|^{2}\lesssim&\frac{1}{\vartheta}\int(\int_{[\frac{s}{\vartheta}]\vartheta}^{s} |\psi_{x}^{n,\mathfrak{s}}
(u,y)|du)^2dy+\int_{[\frac{s}{\vartheta}]\vartheta}^{s}\int |\psi_{x}^{n,\mathfrak{s}}
(t,y)|^2dtdy\\\lesssim&1\wedge (2^{2n}s)\wedge(2^{2n}\vartheta).\endaligned$$
Combining the above estimates we obtain
$$\aligned &\textbf{E}|\langle (\xi_{\varepsilon,\vartheta}-\xi)1_{t\in[0,s]},\psi_{x}^{n,\mathfrak{s}}\rangle|^{2}\lesssim1\wedge (2^{2n}s)\wedge(2^{2n}\varepsilon^2)+1\wedge (2^{2n}s)\wedge(2^{2n}\vartheta).\endaligned$$

Thus, it follows that for $0<\kappa<-\frac{5}{2}-\alpha$,
$$\aligned &\textbf{E}\|(\xi_{\varepsilon,\vartheta}-\xi)1_{t\in[0,s]}\|^{2p}_{\alpha;\mathfrak{R}}\lesssim [\vartheta^{\frac{-\frac{5}{2}-\alpha-\kappa}{2}}+\varepsilon^{-\frac{5}{2}-\alpha-\kappa}]^{2p}s^{\frac{\kappa p}{2}-\frac{5}{2}}.\endaligned$$
Then the required bound (3.4) follows from Kolomogorov's continuity criterion by choosing $p$ large enough.

Now we  prove the second result: it has been obtained in [Hai14, Proposition 9.5] that for every $0<\kappa<-\frac{2\alpha+5}{4}$, the bound
$\textbf{E}\sup_{t\in[0,1]}\|K*\xi_{\varepsilon}(t,\cdot)-K*\xi(t,\cdot)\|_{\alpha+2}\lesssim\varepsilon^{2\kappa}$
holds uniformly over $\varepsilon,\vartheta\in(0,1]$. It suffices to consider $\textbf{E}\sup_{t\in[0,1]}\|K*\xi_{\varepsilon,\vartheta}(t,\cdot)-K*\xi_\varepsilon(t,\cdot)\|_{\alpha+2}$.
We choose the scaling $\bar{\mathfrak{s}}=(1,1,1)$ and
choose a wavelet basis on $\mathbb{R}^3$: $\{\psi_x^{n,\bar{\mathfrak{s}}}=2^{\frac{3}{2}n}\psi(2^n(\cdot-x)),\varphi(\cdot-y), n\geq0, x\in \Lambda^n_{\bar{\mathfrak{s}}}, y\in \Lambda^0_{\bar{\mathfrak{s}}}, \psi\in \bar{\Psi}\}$ as in [Hai14, Section 3.2] on $\mathbb{R}^3$ with $\Lambda_{\bar{\mathfrak{s}}}^n=\{\sum_{j=1}^32^{-n}k_je_j:k_j\in\mathbb{Z}\}$. Set $\bar{\Psi}_\star=\bar{\Psi}\cup\{\varphi\}$. We would like to estimate
$\textbf{E}\|(K*\xi_{\varepsilon,\vartheta}-K*\xi_\varepsilon)(t,\cdot)-(K*\xi_{\varepsilon,\vartheta}-K*\xi_\varepsilon)(s,\cdot)\|^{2p}_{\alpha+2}$ for $t>s\geq0$ and use Kolmogorov's continuity test. Here we only consider the case that $s=0$ for simplicity. For general $s$, we can obtain the desired estimates similarly.
We note that $\|\cdot\|_{\alpha+2}$ on $\mathbb{T}^3$ is equivalent to  the Besov norm $\|\cdot\|_{B^{\alpha+2}_{\infty,\infty}(\mathbb{T}^3)}$, which by [Tri83, Theorem 9.2.1] is equivalent to the weighted Besov norm (cf. [RZZ15, (2.1)]). Moreover, by [Tri06, Theorem 6.15] we have on $\mathbb{T}^3$, for every $p>1$
$$\|f\|_{\alpha+2}^{2p}\lesssim\sum_{\psi\in\bar{\Psi}_\star}
\sum_{n\geq0}\sum_{x\in\Lambda^n_{\bar{\mathfrak{s}}}}2^{2(\alpha+2)pn+|{\bar{\mathfrak{s}}}|np}
|\langle f,\psi_{x}^{n,{\bar{\mathfrak{s}}}}\rangle|^{2p}w(x)^{2p},$$
for $w(x)=(1+|x|^2)^{-2}$,
 which combined with  Gaussian hypercontractivity implies the bound
$$\aligned &\textbf{E}\|(K*\xi_{\varepsilon,\vartheta}-K*\xi_\varepsilon)(t,\cdot)-(K*\xi_{\varepsilon,\vartheta}-K*\xi_\varepsilon)(0,\cdot)\|^{2p}_{\alpha+2}
\\\lesssim&\sum_{\psi\in\bar{\Psi}_\star}
\sum_{n\geq0}\sum_{x\in\Lambda^n_{\bar{\mathfrak{s}}}}2^{2(\alpha+2)pn+|{\bar{\mathfrak{s}}}|np}
(\textbf{E}|\langle (K*\xi_{\varepsilon,\vartheta}-K*\xi_\varepsilon)(t,\cdot)
\\&-(K*\xi_{\varepsilon,\vartheta}-K*\xi_\varepsilon)(0,\cdot),\psi_{x}^{n,{\bar{\mathfrak{s}}}}\rangle|^{2})^pw(x)^{2p}.\endaligned$$
We have the following identity
$$\aligned K*\xi_{\varepsilon,\vartheta}(t,y)=&\sum_{k=-\infty}^{[\frac{t}{\vartheta}]-1}\int_{k\vartheta}^{(k+1)\vartheta}\int K(t-u,y-y_1)\frac{1}{\vartheta}
\int_{k\vartheta}^{(k+1)\vartheta}
\xi_\varepsilon(u_1,y_1)du_1dy_1du\\&+\int_{[\frac{t}{\vartheta}]\vartheta}^{t}\int K(t-u,y-y_1)\frac{1}{\vartheta}
\int_{[\frac{t}{\vartheta}]\vartheta}^{([\frac{t}{\vartheta}]+1)\vartheta}
\xi_\varepsilon(u_1,y_1)du_1dy_1du,\endaligned$$
which implies that
$$\aligned &\textbf{E}|\langle (K*\xi_{\varepsilon,\vartheta}-K*\xi_\varepsilon)(t,\cdot)-(K*\xi_{\varepsilon,\vartheta}-K*\xi_\varepsilon)(0,\cdot),
\psi_{x}^{n,\bar{\mathfrak{s}}}(\cdot)\rangle|^{2}
\\\lesssim&J^1+J^2+J^3+J^4,\endaligned$$
where
$$\aligned J^1:=&\textbf{E}|\langle\sum_{k=0}^{[\frac{t}{\vartheta}]-1}\frac{1}{\vartheta}\int\int_{k\vartheta}^{(k+1)\vartheta}
\int_{k\vartheta}^{(k+1)\vartheta}[K(t-u,\cdot-y_1)-K(t-u_1,\cdot-y_1)]du\\&
\xi_\varepsilon(u_1,y_1)du_1dy_1,\psi_{x}^{n,\bar{\mathfrak{s}}}(\cdot)\rangle|^2,\\J^2:=&\textbf{E}|\langle\frac{1}{\vartheta}\int\int_{[\frac{t}{\vartheta}]\vartheta}^{t} K(t-u,\cdot-y_1)du
\int_{[\frac{t}{\vartheta}]\vartheta}^{([\frac{t}{\vartheta}]+1)\vartheta}
\xi_\varepsilon(u_1,y_1)du_1dy_1,\psi_{x}^{n,\bar{\mathfrak{s}}}(\cdot)\rangle|^2,\\J^3:=&\textbf{E}|\langle \int\int_{[\frac{t}{\vartheta}]\vartheta}^{t}K(t-u_1,\cdot-y_1)
\xi_\varepsilon(u_1,y_1)du_1dy_1,\psi_{x}^{n,\bar{\mathfrak{s}}}(\cdot)\rangle|^2\endaligned$$
and
$$\aligned
J^4:=&\textbf{E}|\langle\int\sum_{k=-\infty}^{-1}\int_{k\vartheta}^{(k+1)\vartheta}\frac{1}{\vartheta}
\int_{k\vartheta}^{(k+1)\vartheta}[K(t-u,\cdot-y_1)-K(t-u_1,\cdot-y_1)-K(-u,\cdot-y_1)\\&+K(-u_1,\cdot-y_1)]du
\xi_\varepsilon(u_1,y_1)du_1dy_1,\psi_{x}^{n,\bar{\mathfrak{s}}}(\cdot)\rangle|^2.\endaligned$$
Now we bound each term separately: For $J^1$ we have
$$\aligned J^1\lesssim&\int\sum_{k=0}^{[\frac{t}{\vartheta}]-1}
\int_{k\vartheta}^{(k+1)\vartheta}\langle\int_{k\vartheta}^{(k+1)\vartheta}\frac{1}{\vartheta}|K(t-u,\cdot-y_1)-K(t-u_1,\cdot-y_1)|du
,|\psi_{x}^{n,\bar{\mathfrak{s}}}(\cdot)|\rangle^2 du_1dy_1,\endaligned$$
We introduce the notation: for $(t,y)\in \mathbb{R}^4$, $\alpha\in\mathbb{R}^+$
\begin{equation}G_0^{(\alpha)}(t,y):=\frac{1}{|t|^{\frac{\alpha}{2}}+|y|^{\alpha}}1_{\{|t|+|y|^2\leq C\}}.\end{equation}
Here $C$ is a constant.
Now we  use [Hai14, Theorem 10.18] to control $|K(t-u,y-y_1)-K(t-u_1,y-y_1)|$ by $\vartheta^{\frac{\delta}{2}} (G_0^{(3+\delta)}(t-u,y-y_1)+G_0^{(3+\delta)}(t-u_1,y-y_1) )$, which implies that
$$\aligned J^1\lesssim&\int\sum_{k=0}^{[\frac{t}{\vartheta}]-1}
\frac{1}{\vartheta^2}\int_{k\vartheta}^{(k+1)\vartheta}\int_{k\vartheta}^{(k+1)\vartheta}\int_{k\vartheta}^{(k+1)\vartheta}\vartheta^\delta (G_0^{(3+\delta)}(t-u,y-y_1)+G_0^{(3+\delta)}(t-u_1,y-y_1) )
\\&(G_0^{(3+\delta)}(t-\tilde{u},\bar{y}-y_1)+G_0^{(3+\delta)}(t-u_1,\bar{y}-y_1) )
|\psi_{x}^{n,\bar{\mathfrak{s}}}(y)\psi_{x}^{n,\bar{\mathfrak{s}}}(\bar{y})|dud\tilde{u}du_1dY.\endaligned$$
Here and in the following we introduce the notation $dY$ to denote  $dyd\bar{y}dy_1$ if there's no confusion. Differently from [Hai14], we calculate the integrals with respect to time and space separately. Observing that each term   on the right hand side of the above inequality only contains at most two of $u, \tilde{u}, u_1$ and using $[|t-u|^{\frac{3}{2}+\frac{\delta}{2}}+|y-y_1|^{3+\delta}]^{-1}\lesssim|t-u|^{-\frac{1-\beta}{2}}|y-y_1|^{-2-\delta-\beta}$ for $0<\beta<1$, we have that $J^1$ can be bounded by
\begin{equation}\aligned&\vartheta^\delta\int
\sum_{k=0}^{[\frac{t}{\vartheta}]-1}
\frac{1}{\vartheta}\int_{k\vartheta}^{(k+1)\vartheta}\int_{k\vartheta}^{(k+1)\vartheta}\big[{|t-u|^{-\frac{1-\beta}{2}}}{|y-y_1|^{-2-\delta-\beta}}
{|t-\tilde{u}|^{-\frac{1-\beta}{2}}}{|\bar{y}-y_1|^{-2-\delta-\beta}}
\\&+{|t-u|^{{-1+\beta}}}{|y-y_1|^{-2-\delta-\beta}}{|\bar{y}-y_1|^{-2-\delta-\beta}}\big]
|\psi_{x}^{n,\bar{\mathfrak{s}}}(y)\psi_{x}^{n,\bar{\mathfrak{s}}}(\bar{y})|dud\tilde{u}1_{\{|y_1|\leq C\}}dY\\
\lesssim&\vartheta^\delta\int\int_{0}^{[\frac{t}{\vartheta}]\vartheta}{|t-u|^{{-1+\beta}}}{|y-y_1|^{-2-\delta-\beta}}{|\bar{y}-y_1|^{-2-\delta-\beta}}
|\psi_{x}^{n,\bar{\mathfrak{s}}}(y)\psi_{x}^{n,\bar{\mathfrak{s}}}(\bar{y})|1_{\{|y_1|\leq C\}}dudY\\\lesssim& \vartheta^\delta|t|^{\beta}\int\int |\psi_{x}^{n,\bar{\mathfrak{s}}}(y)\psi_{x}^{n,\bar{\mathfrak{s}}}(\bar{y})||y-\bar{y}|^{-1-2\delta-2\beta} dyd\bar{y},\endaligned\end{equation}
for $\beta,\delta>0, 2\beta+4\delta<-(2\alpha+5)$. Here in the first inequality we used Young's inequality and in the last inequality we used [Hai14, Lemma 10.14].
For $J^2, J^3$ by similar calculations and the fact that $|K(z)|\lesssim \|z\|_{\mathfrak{s}}^{-3} $  we have
$$\aligned J^2+J^3\lesssim&\int
 \int_{[\frac{t}{\vartheta}]\vartheta}^{t}\int_{[\frac{t}{\vartheta}]\vartheta}^{t}
 G_0^{(3)}(t-u,y-y_1)
G_0^{(3)}(t-\tilde{u},\bar{y}-y_1)dud\tilde{u}\frac{1}{\vartheta}
|\psi_{x}^{n,\bar{\mathfrak{s}}}(y)\psi_{x}^{n,\bar{\mathfrak{s}}}(\bar{y})|dY\\&+\int\int_{[\frac{t}{\vartheta}]\vartheta}^{t}| G_0^{(3)}(t-u_1,y-y_1)
G_0^{(3)}(t-u_1,\bar{y}-y_1)
\psi_{x}^{n,\bar{\mathfrak{s}}}(y)\psi_{x}^{n,\bar{\mathfrak{s}}}(\bar{y})|du_1dY.\endaligned$$
Similar calculations as in (3.6) yield that $J^2+J^3$ can also be bounded by
$$\aligned&\vartheta^\delta|t|^{\beta}\int\int |\psi_{x}^{n,\bar{\mathfrak{s}}}(y)\psi_{x}^{n,\bar{\mathfrak{s}}}(\bar{y})||y-\bar{y}|^{-1-2\delta-2\beta}dyd\bar{y}.\endaligned$$
 Now we consider $J^4$:
$$\aligned J^4\lesssim&
\int\sum_{k=-[\frac{C}{\vartheta}]}^{-1}\frac{1}{\vartheta^2}\int_{k\vartheta}^{(k+1)\vartheta}\int_{k\vartheta}^{(k+1)\vartheta}|K(t-u,y-y_1)-K(t-u_1,y-y_1)
\\&-K(-u,y-y_1)+K(-u_1,y-y_1)|du
\int_{k\vartheta}^{(k+1)\vartheta}|K(t-\tilde{u},\bar{y}-y_1)-K(t-u_1,\bar{y}-y_1)\\&-K(-\tilde{u},\bar{y}-y_1)+K(-u_1,\bar{y}-y_1)|d\tilde{u}
|\psi_{x}^{n,\bar{\mathfrak{s}}}(y)\psi_{x}^{n,\bar{\mathfrak{s}}}(\bar{y})|du_1dY.\endaligned$$
Here we used that $K$ has  compact support. Now we can use [Hai14, Theorem 10.18] to control $|K(t-u,y-y_1)-K(t-u_1,y-y_1)|$ by $\vartheta^{\frac{\delta}{2}} (G^{(3+\delta)}_0(t-u,y-y_1)+G^{(3+\delta)}_0(t-u_1,y-y_1))$ and to control $|K(t-u,y-y_1)-K(-u,y-y_1)|$ by $t^{\frac{\beta}{2}} (G^{(3+\beta)}_0(t-u,y-y_1)+G^{(3+\beta)}_0(-u,y-y_1))$ for $\delta,\beta>0$, which combined  with  interpolation and similar calculations as in (3.6) imply that
$$\aligned J^4\lesssim&\vartheta^\delta|t|^{\beta}\int\int |\psi_{x}^{n,\bar{\mathfrak{s}}}(y)\psi_{x}^{n,\bar{\mathfrak{s}}}(\bar{y})||y-\bar{y}|^{-1-4\delta-2\beta} dyd\bar{y}.\endaligned$$
Combining  the above estimates  we obtain that
$$\aligned &\textbf{E}|\langle (K*\xi_{\varepsilon,\vartheta}-K*\xi_\varepsilon)(t,\cdot)-(K*\xi_{\varepsilon,\vartheta}-K*\xi_\varepsilon)(0,\cdot),\psi_{x}^{n,\bar{\mathfrak{s}}}\rangle|^{2}
\\\lesssim& \vartheta^\delta|t|^{\beta}\int\int |\psi_{x}^{n,\bar{\mathfrak{s}}}(y)\psi_{x}^{n,\bar{\mathfrak{s}}}(\bar{y})||y-\bar{y}|^{-1-4\delta-2\beta} dyd\bar{y}\\\lesssim& \vartheta^\delta|t|^{\beta}2^{-3n+n(1+4\delta+2\beta)},\endaligned$$
where $\beta,\delta>0, 2\beta+4\delta<-(2\alpha+5)$.
Thus, the above estimates yield that $$\aligned &\textbf{E}\|(K*\xi_{\varepsilon,\vartheta}-K*\xi_\varepsilon)(t,\cdot)-(K*\xi_{\varepsilon,\vartheta}-K*\xi_\varepsilon)(0,\cdot)\|^{2p}_{\alpha+2}\\\lesssim&\sum_{\psi\in\bar{\Psi}_\star}
\sum_{n\geq0}2^{2(\alpha+\frac{5}{2})np+|\bar{\mathfrak{s}}|np+|\bar{\mathfrak{s}}|n}
\vartheta^{p\delta}|t|^{\beta p}2^{-3np+np(4\delta+2\beta)},\endaligned$$
and the results follow from Kolmogorov's continuity test (in time) if we choose $p$ sufficiently large.
$\hfill\Box$
\vskip.10in
In [Hai14, Theorem 10.22] a random model $\hat{Z}$ has been obtained by taking the limit of the models associated with the convolution approximation $\xi_\varepsilon$. Define $\Phi=\mathcal{R}\mathcal{S}^L(u_0,\hat{Z})$. Then $\Phi$ is the local solution to the dynamical $\Phi^4_3$ model.
In our case we also  have the following main convergence result at the level of models:
\vskip.10in
\th{Theorem 3.7} Let $\mathfrak{T}_F$ be the regularity structure associated to the dynamical $\Phi_3^4$ model, and $\rho(t,x)=\rho_1(t)\rho_2(x)$, let $\xi_{\varepsilon,\vartheta}$ be as in Subsection 2.2  and let $Z_{\varepsilon,\vartheta}$ be the associated model. Then there exist choices of $C_1^{(\varepsilon,\vartheta)}(t), C_2^{(\varepsilon,\vartheta)}(t)$ such that $\hat{Z}_{\varepsilon,\vartheta}=Z_{\varepsilon,\vartheta}^M\rightarrow\hat{Z}$ in probability.

More precisely, for any $\kappa<-\frac{5}{2}-\alpha$, any compact set $\mathfrak{R}$, and any $\gamma<r$ one has that the bound
$$\textbf{E}|||Z_{\varepsilon,\vartheta}^M;\hat{Z}|||_{\gamma;\mathfrak{R}}\lesssim \varepsilon^\kappa+\vartheta^{\kappa/2},$$
holds uniformly over $\varepsilon, \vartheta\in (0,1]$.
\vskip.10in
The proof of Theorem 3.7 is the content of Section 4 and the Appendix below.

\vskip.10in
\th{Remark 3.8} If $\rho(t,x)=\delta(t)\rho_2(x)$ for the Dirac distribution $\delta$, the convergence results in Proposition 3.6 and Theorem 3.7 still hold. In fact, if $|K(z)|\lesssim \|z\|_\mathfrak{s}^\zeta$ for $-4<\zeta<0$ it is sufficient to control $|K*\rho_{2,\varepsilon}-K|$ by $(t^{-\frac{\delta}{2}}\varepsilon^{\zeta-\bar{\zeta}}|x|^{\bar{\zeta}+\delta})\wedge|t|^{\frac{\zeta}{2}}$ for $\bar{\zeta}+\delta>-3$. By this  and similar calculations as in Section 4 we could also deduce the results. Here $\rho_{2,\varepsilon}(y)=\varepsilon^{-3}\rho_2(\frac{y}{\varepsilon})$.

\vskip.10in

We now have all the tools in place to prove the main convergence result of
this article.
\vskip.10in

\no \emph{Proof of Theorem 1.1} The proof of the theorem is essentially a collection of the results of this paper.
As obtained in Proposition 3.5, $\mathcal{R}\mathcal{S}^L(u_0,Z_{\varepsilon,\vartheta}^M)=\Phi_{\varepsilon,\vartheta}$. Define $\Phi=\mathcal{R}\mathcal{S}^L(u_0,\hat{Z})$. By the continuity of the map $\mathcal{R}$ and Proposition 2.8, Theorem 3.7, we obtain that
there exists a sequence of random times
$\tau_L$ converging to the explosion time $\tau$ of $\Phi$ such that
$$\sup_{t\in[0,\tau_L]}\|\Phi_{\varepsilon,\vartheta}-\Phi\|_\eta\rightarrow^P0, \textrm{ as }\varepsilon, \vartheta\rightarrow0.$$

\section{Convergence of the renormalised model}
 In the previous section we have defined the renormalised models $\hat{Z}_{\varepsilon,\vartheta}=(\hat{\Pi}^{(\varepsilon,\vartheta)}, \hat{\Gamma}^{(\varepsilon,\vartheta)})$. The goal of this section is to obtain the convergence of the renormalised models. The proof follows by a similar argument as in the proof of [Hai14, Theorem 10.22], if we can prove the following lemmas. For the completeness of the paper we put the proof of Theorem 3.7 in the Appendix. In the following we prove these lemmas.
First, we introduce  the following notations:

Define for $(t,y), (t_2,y_2)\in\mathbb{R}^4$ \begin{equation}\aligned K_{\varepsilon,\vartheta}(t,y,t_2,y_2):=&\sum_{k=-\infty}^{[\frac{t}{\vartheta}]-1}\int_{k\vartheta}^{(k+1)\vartheta}
\int K(t-u,y-y_1)\frac{1}{\vartheta}
\int_{k\vartheta}^{(k+1)\vartheta}
\rho_\varepsilon(u_1-t_2,y_1-y_2)du_1dy_1du\\&+\int_{[\frac{t}{\vartheta}]\vartheta}^{t}\int K(t-u,y-y_1)\frac{1}{\vartheta}
\int_{[\frac{t}{\vartheta}]\vartheta}^{([\frac{t}{\vartheta}]+1)\vartheta}
\rho_\varepsilon(u_1-t_2,y_1-y_2)du_1dy_1du\\:=&K_{\varepsilon,\vartheta}^{(1)}(t,y,t_2,y_2)+K_{\varepsilon,\vartheta}^{(2)}(t,y,t_2,y_2),\endaligned\end{equation} \begin{equation}\aligned K_\varepsilon(t-t_2,y-y_2):=&K*\rho_\varepsilon=\int\int_{-\infty}^{[\frac{t}{\vartheta}]\vartheta}K(t-u,y-y_1)\rho_{\varepsilon}(u-t_2,y_1-y_2)dudy_1
\\&+\int\int_{[\frac{t}{\vartheta}]\vartheta}^tK(t-u,y-y_1)\rho_{\varepsilon}(u-t_2,y_1-y_2)dudy_1\\:=&K_\varepsilon^{(1)}(t,y,t_2,y_2)+K_\varepsilon^{(2)}
(t,y,t_2,y_2).\endaligned\end{equation}
Then for $\tau=\mathcal{I}(\Xi)=\Psi$ we have
$$\aligned ({\hat{\Pi}}^{(\varepsilon,\vartheta)}\Psi)(z)=&K*\xi_{\varepsilon,\vartheta}(z)=\int K_{\varepsilon,\vartheta}(z,z_1)\xi(z_1)dz_1.\endaligned$$
For $\varphi$ smooth and $x\in\mathbb{R}^4$ we have that
$$\aligned &E|\langle K*\xi_{\varepsilon,\vartheta},\varphi_{x}^\lambda\rangle|^{2}=\int\int f^{(\varepsilon,\vartheta)}(z,\bar{z})
\varphi_{x}^\lambda(z)\varphi_{x}^\lambda(\bar{z})dzd\bar{z}.\endaligned$$
Here for $z=(t,y), \bar{z}=(\bar{t},\bar{y})$
$$\aligned &f^{(\varepsilon,\vartheta)}(z,\bar{z}):=\sum_{i,j=1}^2J^{ij}(z,\bar{z}), \endaligned$$
with
$$\aligned
J^{ij}(z,\bar{z})=&\int K_{\varepsilon,\vartheta}^{(i)}(z,z_1)K_{\varepsilon,\vartheta}^{(j)}(\bar{z},z_1)dz_1, i,j=1,2.\endaligned$$

In the following we prove estimates for $f^{(\varepsilon,\vartheta)}$ and $K_{\varepsilon,\vartheta}-K_\varepsilon$. Recall that $\rho(t,x)=\rho_1(t)\rho_2(x)$. We  first give an estimate for the convolution $K*\rho_{2,\varepsilon}$ with respect to space, which is required for the estimate of $f^{(\varepsilon,\vartheta)}$. Here $\rho_{2,\varepsilon}(y)=\varepsilon^{-3}\rho_2(\frac{y}{\varepsilon})$. By a similar argument as the proof in [Hai14, Lemma 10.17] we obtain:
\vskip.10in

\th{Lemma 4.1} If $|K(z)|\lesssim \|z\|_{\mathfrak{s}}^{\zeta}$ for $\zeta\in (-4,0)$, then
$$|K*\rho_{2,\varepsilon}(z)|\leq Ct^{-\frac{\delta}{2}}(\|z\|_{\mathfrak{s}}^{\zeta+\delta}\wedge\varepsilon^{\zeta+\delta}),$$
for $ 0<\delta<1, 0>\zeta+\delta>-3$.

\proof We can write $$K*\rho_{2,\varepsilon}(t,x)=\int K(t,x-y)\rho_{2,\varepsilon}(y)dy.$$
We use the notation $z=(t,x)$. $|K(t,x-y)|$ can be bounded by $C|t|^{\frac{\zeta}{2}}$, and
 $$|K*\rho_{2,\varepsilon}(t,x)|\leq C|t|^{\frac{\zeta}{2}}$$
 follows from the fact that $\rho_{2,\varepsilon}$ integrates to $1$.
Without loss of generality we assume that $\rho_{2}$ is supported in the set $\{x:|x|\leq 1\}$. For $|x|\geq 2\varepsilon$, we have $|x-y|\geq \frac{|x|}{2}$, which implies that for $0<\delta<1, 0>\zeta+\delta>-3$,
 $$|K*\rho_{2,\varepsilon}(t,x)|\leq Ct^{-\frac{\delta}{2}}|x|^{{\zeta+\delta}} \leq Ct^{-\frac{\delta}{2}}\varepsilon^{\zeta+\delta}.$$
 For $|x|\leq 2\varepsilon$ we use the fact that $|\rho_{2,\varepsilon}|$ is bounded by a constant multiple of $\varepsilon^{-3}$
 $$|K*\rho_{2,\varepsilon}(t,x)|\lesssim \varepsilon^{-3}\int_{|y|\leq 3\varepsilon} t^{-\frac{\delta}{2}}|y|^{\zeta+\delta} dy\lesssim Ct^{-\frac{\delta}{2}}\varepsilon^{{\zeta+\delta}}\lesssim Ct^{-\frac{\delta}{2}}|x|^{{\zeta+\delta}}.$$
 Combining all the estimates the result follows.$\hfill\Box$

\vskip.10in
In the following we prove some  useful estimates for $f^{(\varepsilon,\vartheta)}$ and  $K_{\varepsilon,\vartheta}^{(i)}-K_{\varepsilon}^{(i)}$, which are used  in the proof of Theorem 3.7.

\vskip.10in
\th{Lemma 4.2}

(i) For every $\delta>0$
$$|f^{(\varepsilon,\vartheta)}(z,\bar{z})|\lesssim\|z-\bar{z}\|_\mathfrak{s}^{-1-\delta}.$$

(ii)$$\aligned &|f^{(\varepsilon,\vartheta)}(z,\bar{z})-f(z-\bar{z})|
\lesssim(\vartheta^{\kappa}+\varepsilon^{2\kappa})\|z-\bar{z}\|_\mathfrak{s}^{-1-{2\kappa}-\delta}\endaligned$$
holds uniformly over $\varepsilon,\vartheta\in(0,1]$, provided that $\kappa<1$ and that $\delta>0$,  where $f(z-\bar{z})=K*K(z-\bar{z})$.

(iii)  For $i=1,2$
$$\aligned &|\langle (K_{\varepsilon,\vartheta}^{(i)}-K_\varepsilon^{(i)})(z,\cdot),(K_{\varepsilon,\vartheta}^{(i)}-K_\varepsilon^{(i)})(\bar{z},\cdot)\rangle|
\lesssim\vartheta^{\kappa}\|z-\bar{z}\|_\mathfrak{s}^{-1-{2\kappa}-\delta},\endaligned$$
holds uniformly over $\varepsilon,\vartheta\in(0,1]$, provided that $\kappa<1$ and that $\delta>0$.

\proof In the following we use the notations $z=(t,y), \bar{z}=(\bar{t},\bar{y})$. Consider (i) first: We consider the integral w.r.t. space and time separately: Since $\rho_1$ has compact support,   there exists some constant $C_0$ such that for $\rho_{1,\varepsilon}(t):=\varepsilon^{-2}\rho_1(\frac{t}{\varepsilon^2})$,
$\int_{k\vartheta}^{(k+1)\vartheta}\int_{k_1\vartheta}^{(k_1+1)\vartheta}\rho_{1,\varepsilon}*\rho_{1,\varepsilon}(u_1-u_2)du_1du_2\neq0$
if and only if $|k-k_1|\leq \frac{C_0\varepsilon^2}{\vartheta}+1$, and in this case
\begin{equation}|\int_{k\vartheta}^{(k+1)\vartheta}\int_{k_1\vartheta}^{(k_1+1)\vartheta}\rho_{1,\varepsilon}*\rho_{1,\varepsilon}(u_1-u_2)du_1du_2|\lesssim\vartheta(1\wedge
\frac{\vartheta}{\varepsilon^2}),\end{equation}
which implies that
\begin{equation}\aligned J^{11}\lesssim&\sum_{k=-[\frac{C}{\vartheta}]}^{[\frac{t}{\vartheta}]-1}\sum_{k_1=-[\frac{C}{\vartheta}], |k_1-k|\leq \frac{C_0\varepsilon^2}{\vartheta}+1}^{[\frac{\bar{t}}{\vartheta}]-1}\int\frac{1}{\vartheta}\int_{k\vartheta}^{(k+1)\vartheta}|t-u|^{-\frac{\delta}{2}}
[G_0^{(3-\delta)}(t-u,y-y_1)\wedge\varepsilon^{-3+\delta}]du
\\&\int_{k_1\vartheta}^{(k_1+1)\vartheta}|\bar{t}-\tilde{u}|^{-\frac{\delta}{2}}[G_0^{(3-\delta)}(\bar{t}-\tilde{u},\bar{y}-y_1)\wedge\varepsilon^{-3+\delta}]
d\tilde{u}dy_1(1\wedge
\frac{\vartheta}{\varepsilon^2}),\endaligned\end{equation}
for $\delta>0$, where we used  Lemma 4.1 and $G_0^{(3-\delta)}$ is defined as in (3.5). Now we consider this term in the following three cases:

Case I: $\bar{t}-t\geq 2C_0\varepsilon^2+4\vartheta$. Since $|u-{\tilde{u}}|\leq C_0\varepsilon^2+2\vartheta$, we deduce that
\begin{equation}\bar{t}-\tilde{u}>t-u.\end{equation}
Furthermore, we  have that $$\bar{t}-\tilde{u}=\bar{t}-t+t-u+u-\tilde{u}\geq C_0\varepsilon^2+2\vartheta\geq \tilde{u}-u+u-t=\tilde{u}-t,$$ which  implies that
\begin{equation}|\bar{t}-\tilde{u}|\geq\frac{|\bar{t}-t|}{2}.\end{equation}
By (4.4), (4.5), (4.6) and using \begin{equation}[|t-u|^{\frac{3-\delta}{2}}+|y-y_1|^{3-\delta}]^{-1}\lesssim (t-u)^{-a}|y-y_1|^{-b}\end{equation} for $a,b>0$ satisfying  $2a+b=3-\delta$,  we obtain that
\begin{equation}\aligned J^{11}\lesssim&\sum_{k=-[\frac{C}{\vartheta}]}^{[\frac{t}{\vartheta}]-1}\int\int_{k\vartheta}^{(k+1)\vartheta}
\bigg({(t-u)^{-1+\delta/4}|y-y_1|^{-1-\delta/2}}
{(\bar{t}-t)^{-\frac{1}{2}-\frac{\delta}{2}}|\bar{y}-y_1|^{-2+\delta}}\bigg)\\&\bigwedge \bigg({(t-u)^{-1+\delta/2}|y-y_1|^{-1-2\delta}}
{|\bar{y}-y_1|^{-3+\delta}}\bigg)du1_{|y_1|\leq C}dy_1\\\lesssim&|t-\bar{t}|^{-\frac{1}{2}-\frac{\delta}{2}}\wedge|y-\bar{y}|^{-1-\delta},\endaligned\end{equation}
is valid for every $\delta>0$, where we first used (4.4), (4.7) and then (4.5), (4.6) in the first inequality and used [Hai14, Lemma 10.14] in the last inequality.

Case II: $t-\bar{t}\geq 2C_0\varepsilon^2+4\vartheta$. Similarly as Case I.

Case III: $|\bar{t}-t|\leq 2C_0\varepsilon^2+4\vartheta$. We have that for every $\delta>0$
$$\frac{1}{\vartheta}\int_{k_1\vartheta}^{(k_1+1)\vartheta}{(\bar{t}-\tilde{u})^{-\frac{1}{2}-\frac{\delta}{2}}}d\tilde{u}\lesssim \vartheta^{-\frac{1}{2}-\frac{\delta}{2}},$$
which combined with (4.7) implies that
\begin{equation}\aligned J^{11}\lesssim&\sum_{k=-[\frac{C}{\vartheta}]}^{[\frac{t}{\vartheta}]-1}\int\int_{k\vartheta}^{(k+1)\vartheta}{(t-u)^{-1+\delta/4}|y-y_1|^{-1-\delta/2}}
{\vartheta^{-\frac{1}{2}-\frac{\delta}{2}}|\bar{y}-y_1|^{-2+\delta}}du1_{\{|y_1|\leq C\}}dy_1\\\lesssim&\vartheta^{-\frac{1}{2}-\frac{\delta}{2}}.
\endaligned\end{equation}
holds for every $\delta>0$, where we used [Hai14, Lemma 10.14] in the last inequality.  (4.4) and (4.7) also imply that
 \begin{equation}\aligned J^{11}\lesssim&(\frac{\vartheta}{\varepsilon^2}\wedge1)\sum_{k=-[\frac{C}{\vartheta}]}^{[\frac{t}{\vartheta}]-1}
\sum_{k_1=-[\frac{C}{\vartheta}], |k_1-k|\leq C_0\frac{\varepsilon^2}{\vartheta}+1}^{[\frac{\bar{t}}{\vartheta}]-1}\int\frac{1}{\vartheta}\int_{k\vartheta}^{(k+1)\vartheta}
\int_{k_1\vartheta}^{(k_1+1)\vartheta} [{(t-u)^{-\frac{1}{2}+\frac{\delta}{4}}|y-y_1|^{-2-\frac{\delta}{2}}}
\\&{(\bar{t}-\tilde{u})^{-\frac{1}{2}+\frac{\delta}{4}}}{|\bar{y}-y_1|^{-2-\frac{\delta}{2}}}]dud\tilde{u}1_{\{|y_1|\leq C\}}dy_1\lesssim|y-\bar{y}|^{-1-\delta},\endaligned\end{equation}
 where we used [Hai14, Lemma 10.14] and ${(t-u)^{-\frac{1}{2}+\frac{\delta}{4}}}{(\bar{t}-\tilde{u})^{-\frac{1}{2}+\frac{\delta}{4}}}\leq {(t-u)^{-1+\frac{\delta}{2}}}+{(\bar{t}-\tilde{u})^{-1+\frac{\delta}{2}}}$ in the last inequality. Moreover,  by interpolation we have that
 $$|t-u|^{-\frac{\delta}{2}}
[G_0^{(3-\delta)}(t-u,y-y_1)\wedge\varepsilon^{-3+\delta}]\lesssim{(t-u)^{-\frac{1}{2}+\frac{\delta}{8}}
|y-y_1|^{-\frac{3}{2}+\frac{\delta}{4}}}\varepsilon^{-\frac{1}{2}-\frac{\delta}{2}},$$
  which combined with (4.4)  yields that
 \begin{equation}\aligned J^{11}\lesssim&(\frac{\vartheta}{\varepsilon^2}\wedge1)\sum_{k=-[\frac{C}{\vartheta}]}^{[\frac{t}{\vartheta}]-1}
\sum_{k_1=-[\frac{C}{\vartheta}], |k_1-k|\leq C_0\frac{\varepsilon^2}{\vartheta}+1}^{[\frac{\bar{t}}{\vartheta}]-1}\int\frac{1}{\vartheta}\int_{k\vartheta}^{(k+1)\vartheta}
\int_{k_1\vartheta}^{(k_1+1)\vartheta} [{(t-u)^{-\frac{1}{2}+\frac{\delta}{8}}|y-y_1|^{-\frac{3}{2}+\frac{\delta}{4}}}\varepsilon^{-1-\delta}
\\&{(\bar{t}-\tilde{u})^{-\frac{1}{2}+\frac{\delta}{8}}}{|\bar{y}-y_1|^{-\frac{3}{2}+\frac{\delta}{4}}}]dud\tilde{u}1_{\{|y_1|\leq C\}}dy_1\lesssim\varepsilon^{-1-\delta},\endaligned\end{equation}
 where we used [Hai14, Lemma 10.14] and ${(t-u)^{-\frac{1}{2}+\frac{\delta}{8}}}{(\bar{t}-\tilde{u})^{-\frac{1}{2}+\frac{\delta}{8}}}\leq {(t-u)^{-1+\frac{\delta}{4}}}+{(\bar{t}-\tilde{u})^{-1+\frac{\delta}{4}}}$ in the last inequality.
Combining (4.9)-(4.11) we obtain that $$\aligned J^{11}\lesssim\varepsilon^{-1-\delta}\wedge\vartheta^{-\frac{1}{2}-\frac{\delta}{2}}\wedge|y-\bar{y}|^{-1-\delta}\lesssim|t-\bar{t}|^{-\frac{1}{2}-\frac{\delta}{2}}\wedge|y-\bar{y}|^{-1-\delta}\lesssim&\|z-\bar{z}\|^{-1-\delta}_{\mathfrak{s}},
\endaligned$$
is valid for every $\delta>0$.

We now turn to $J^{22}$. $J^{22}\neq 0$ if and only if $|t-\bar{t}|\leq 2C_0\varepsilon^2+4\vartheta$. Lemma 4.1, (4.3) and similar arguments as in (4.9-4.11) imply that for every $\delta>0$
\begin{equation}\aligned J^{22}\lesssim&\bigg(\int\int_{[\frac{t}{\vartheta}]\vartheta}^{t}{(t-u)^{-1+\delta/4}|y-y_1|^{-1-\frac{\delta}{2}}}
{\vartheta^{-\frac{1+\delta}{2}}|\bar{y}-y_1|^{-2+\delta}}du1_{\{|y_1|\leq C\}}dy_1\bigg)\\&\bigwedge\bigg(\int\frac{1}{\vartheta}\int_{[\frac{t}{\vartheta}]\vartheta}^{t}\int_{[\frac{\bar{t}}{\vartheta}]\vartheta}^{\bar{t}}
[{(t-u)^{-\frac{1}{2}+\frac{\delta}{4}}
|y-y_1|^{-2-\frac{\delta}{2}}}
{(\bar{t}-\tilde{u})^{-\frac{1}{2}+\frac{\delta}{4}}|\bar{y}-y_1|^{-2-\frac{\delta}{2}}}]\\&\wedge [{(t-u)^{-\frac{1}{2}+\frac{\delta}{8}}|y-y_1|^{-\frac{3}{2}+\frac{\delta}{4}}}\varepsilon^{-1-\delta}
{(\bar{t}-\tilde{u})^{-\frac{1}{2}+\frac{\delta}{8}}}{|\bar{y}-y_1|^{-\frac{3}{2}+\frac{\delta}{4}}}]d\tilde{u}du1_{\{|y_1|\leq C\}}dy_1\bigg) \\\lesssim&|t-\bar{t}|^{-\frac{1}{2}-\frac{\delta}{2}}\wedge|y-\bar{y}|^{-1-\delta}
.\endaligned\end{equation}
$J^{12}, J^{21}$ can be estimated similarly.
Thus (i) follows.

(ii) We have
$$\aligned &|f^{(\varepsilon,\vartheta)}(z,\bar{z})-f(z-\bar{z})|\\\lesssim&|f^{(\varepsilon,\vartheta)}(z,\bar{z})-f^{(\varepsilon)}(z-\bar{z})|
+|f^{(\varepsilon)}(z-\bar{z})-f(z-\bar{z})|
\endaligned$$
where $f^{(\varepsilon)}(z-\bar{z})=K_\varepsilon*K_\varepsilon(z-\bar{z}), K_\varepsilon=K*\rho_\varepsilon $.
Here $|f^{(\varepsilon,\vartheta)}(z,\bar{z})-f^{(\varepsilon)}(z-\bar{z})|$ can be separated as $J^{ij}_1, i,j=1,2,$ with
$$\aligned
J^{ij}_1=&\bigg|\int (K_{\varepsilon,\vartheta}^{(i)}(z,z_1)K_{\varepsilon,\vartheta}^{(j)}(\bar{z},z_1)-K_{\varepsilon}^{(i)}(z,z_1)K_{\varepsilon}^{(j)}(\bar{z},z_1))dz_1\bigg|
,\endaligned$$
where $K_{\varepsilon,\vartheta}^{(i)}$ and $K_{\varepsilon}^{(i)}$ are defined as in (4.1) and (4.2).
Each term can be estimated as in the proof of (i).  We take $J^{11}$ as an example:
$$\aligned J^{11}_1=&\bigg|\sum_{k=-\infty}^{[\frac{t}{\vartheta}]-1}\sum_{k_1=-\infty}^{[\frac{\bar{t}}{\vartheta}]-1}\frac{1}{\vartheta^2}\int\int\int_{k\vartheta}^{(k+1)\vartheta}
\int_{k_1\vartheta}^{(k_1+1)\vartheta}
\int_{k\vartheta}^{(k+1)\vartheta}\int_{k_1\vartheta}^{(k_1+1)\vartheta}\big[(K(t-u,y-y_1)\\&-K(t-u_1,y-y_1))
K(\bar{t}-\tilde{u},\bar{y}-y_2)+K(t-u_1,y-y_1)
(K(\bar{t}-\tilde{u},\bar{y}-y_2)\\&-K(\bar{t}-u_2,\bar{y}-y_2))\big]dud\tilde{u}
\rho_\varepsilon*\rho_\varepsilon(u_1-u_2,y_1-y_2)du_1du_2dy_1dy_2\bigg|.\endaligned$$
By [Hai14, Lemma 10.18] we can control
$|K(t-u,y-y_1)-K(t-u_1,y-y_1)|$ by $\vartheta^\kappa (G^{(3+2\kappa)}_0(t-u,y-y_1)+G^{(3+2\kappa)}_0(t-u_1,y-y_1))$ provided that $\kappa>0$. Observing that each term   on the right hand side of the above inequality only contains at most two of $u, \tilde{u}, u_1, u_2$ and using (4.3) as well as  Lemma 4.1 we obtain that $J^{11}_1$ can be bounded by a term similar as in (4.4). Then by similar arguments as the estimates for (4.4), we deduce  that $J^{11}_1\lesssim\vartheta^{\kappa}\|z-\bar{z}\|_\mathfrak{s}^{-1-{2\kappa}-\delta}$. The other terms can be estimated similarly, which implies that
$$|f^{(\varepsilon,\vartheta)}(z,\bar{z})-f^{(\varepsilon)}(z-\bar{z})|
\lesssim\vartheta^{\kappa}\|z-\bar{z}\|_\mathfrak{s}^{-1-{2\kappa}-\delta},$$
holds uniformly over $\varepsilon,\vartheta\in(0,1]$, provided that $\kappa<1$ and that $\delta>0$. Since $K$ and $K_\varepsilon$ are of order $-3$, by [Hai14, Lemma 10.17] we obtain that
$$|f^{(\varepsilon)}(z-\bar{z})-f(z-\bar{z})|\lesssim\varepsilon^{2\kappa}\|z-\bar{z}\|_\mathfrak{s}^{-1-{2\kappa}-\delta}.$$
Combining the above estimates we deduce (ii)  easily.

(iii) We have for $z=(t,y), \bar{z}=(\bar{t},\bar{y})$
$$\aligned &|\langle (K_{\varepsilon,\vartheta}^{(1)}-K_\varepsilon^{(1)})(z,\cdot),(K_{\varepsilon,\vartheta}^{(1)}-K_\varepsilon^{(1)})(\bar{z},\cdot)\rangle|
\\=&\bigg|\sum_{k=-\infty}^{[\frac{t}{\vartheta}]-1}\sum_{k_1=-\infty}^{[\frac{\bar{t}}{\vartheta}]-1}\int\int\int_{k\vartheta}^{(k+1)\vartheta}
\int_{k_1\vartheta}^{(k_1+1)\vartheta}\int_{k\vartheta}^{(k+1)\vartheta}\int_{k_1\vartheta}^{(k_1+1)\vartheta}(K(t-u,y-y_1)-K(t-u_1,y-y_1))\\&
(K(\bar{t}-\tilde{u},\bar{y}-y_2)-K(\bar{t}-u_2,\bar{y}-y_2))dud\tilde{u}\frac{1}{\vartheta^2}
\rho_\varepsilon*\rho_\varepsilon(u_1-u_2,y_1-y_2)du_1du_2dy_1dy_2\bigg|.\endaligned$$
By [Hai14, Lemma 10.18] we can control
$|K(t-u,y-y_1)-K(t-u_1,y-y_1)|$ by $\vartheta^\kappa (G^{(3+2\kappa)}_0(t-u,y-y_1)+G^{(3+2\kappa)}_0(t-u_1,y-y_1))$.
Then by similar arguments as in the proof of (ii) we obtain that $$|\langle (K_{\varepsilon,\vartheta}^{(1)}-K_\varepsilon^{(1)})(z,\cdot),(K_{\varepsilon,\vartheta}^{(1)}-K_\varepsilon^{(1)})(\bar{z},\cdot)
\rangle|\lesssim\vartheta^{\kappa}\|z-\bar{z}\|_\mathfrak{s}^{-1-{2\kappa}-\delta}.$$
For $i=2$ we can argue similarly as for estimating $J^{22}$ in (i) and the result follows.
$\hfill\Box$
\vskip.10in
\th{Lemma 4.3} The following holds:

(i) For every $\delta>0$
$$\aligned &|f^{(\varepsilon,\vartheta)}(z,\bar{z}_1)-f^{(\varepsilon,\vartheta)}(z,\bar{z}_2)|\lesssim\|\bar{z}_1-\bar{z}_2\|_\mathfrak{s}^\delta
(\|z-\bar{z}_1\|_\mathfrak{s}^{-1-2\delta}+\|z-\bar{z}_2\|_\mathfrak{s}^{-1-2\delta}).\endaligned$$

(ii) For $i=1, 2$, $$\aligned&|\langle (K_{\varepsilon,\vartheta}^{(i)}-K_\varepsilon^{(i)})(z,\cdot),(K_{\varepsilon,\vartheta}^{(i)}-K_\varepsilon^{(i)})(\bar{z}_1,\cdot)
-(K^{(i)}_{\varepsilon,\vartheta}-K^{(i)}_\varepsilon)(\bar{z}_2,\cdot)\rangle|\\\lesssim& \vartheta^\kappa\|\bar{z}_1-\bar{z}_2\|_\mathfrak{s}^\delta(\|\bar{z}_1-z\|_\mathfrak{s}^{-1-2\delta-2\kappa}
 +\|\bar{z}_2-z\|_\mathfrak{s}^{-1-2\delta-2\kappa})\endaligned$$
 holds uniformly over $\varepsilon,\vartheta\in(0,1]$, provided that $\kappa<1$ and that $\delta>0$.

\proof Without loss of generality for $\bar{z}_1=(\bar{t}_1,\bar{y}_1), \bar{z}_2=(\bar{t}_2,\bar{y}_2)$ we suppose that $\bar{t}_1\leq \bar{t}_2$.
$|f^{(\varepsilon,\vartheta)}(z,\bar{z}_1)-f^{(\varepsilon,\vartheta)}(z,\bar{z}_2)|$ can be separated into two kinds of terms: one is similar to $J^{ij}$ in the proof of Lemma 4.2 with one of $K(\bar{t},\bar{y})$ in $J^{ij}$ replaced by $K(\bar{t}_1,\bar{y})-K(\bar{t}_2,\bar{y})$; the other is the corresponding terms from $\bar{t}_1$ to $\bar{t}_2$. For the first case, we apply [Hai14, Lemma 10.18] to deduce that $|K(\bar{t}_1,\bar{y})-K(\bar{t}_2,\bar{y})|\lesssim |\bar{t}_1-\bar{t}_2|^\frac{\delta}{2}( G_0^{(3+\delta)}(\bar{t}_1,\bar{y})+G_0^{(3+\delta)}(\bar{t}_2,\bar{y}))$, which combined with  similar arguments as in the proof of Lemma 4.2 implies the desired estimates.
Now we give the calculations for the most complicated term in the second case and the other terms can be handled similarly. Define
$$\aligned J=&\bigg|\sum_{k=-\infty}^{[\frac{t}{\vartheta}]-1}\sum_{k_1=[\frac{\bar{t}_1}{\vartheta}]+1}^{[\frac{\bar{t}_2}{\vartheta}]-1}\int\int\int_{k\vartheta}^{(k+1)\vartheta}K(t-u,y-y_1)du
\int_{k_1\vartheta}^{(k_1+1)\vartheta}K(\bar{t}_2-\tilde{u},\bar{y}_2-y_2)d\tilde{u}\\&\frac{1}{\vartheta^2}
\int_{k\vartheta}^{(k+1)\vartheta}\int_{k_1\vartheta}^{(k_1+1)\vartheta}
\rho_\varepsilon*\rho_\varepsilon(u_1-u_2,y_1-y_2)du_1du_2dy_1dy_2\bigg|.\endaligned$$
By (4.3) and Lemma 4.1, similarly as in (4.4) we have that for $\delta>0$
\begin{equation}\aligned J\lesssim&\sum_{k=-[\frac{C}{\vartheta}]}^{[\frac{t}{\vartheta}]-1}\sum_{k_1=[\frac{\bar{t}_1}{\vartheta}]+1,|k_1-k|\leq C_0\frac{\varepsilon^2}{\vartheta}+1}^{[\frac{\bar{t}_2}{\vartheta}]-1}
\int\int_{k\vartheta}^{(k+1)\vartheta}|t-u|^{-\delta/2}[G^{(3-\delta)}_0(t-u,y-y_1)\wedge\varepsilon^{-3+\delta}]du
\\&\int_{k_1\vartheta}^{(k_1+1)\vartheta}|\bar{t}_2-\tilde{u}|^{-\delta/2}[G^{(3-\delta)}_0(\bar{t}_2-\tilde{u},\bar{y}_2-y_1)\wedge\varepsilon^{-3+\delta}]d\tilde{u}\frac{1}{\vartheta}
dy_1(\frac{\vartheta}{\varepsilon^2}\wedge1).\endaligned\end{equation}
Then by  (4.7) and
$|\bar{t}_2-\bar{t}_1|\geq |\bar{t}_2-\tilde{u}|$ we obtain
$$\aligned J\lesssim&\sum_{k=-[\frac{C}{\vartheta}]}^{[\frac{t}{\vartheta}]-1}\sum_{k_1=[\frac{\bar{t}_1}{\vartheta}]+1, |k_1-k|\leq C_0\frac{\varepsilon^2}{\vartheta}}^{[\frac{\bar{t}_2}{\vartheta}]-1}\int\int_{k\vartheta}^{(k+1)\vartheta}\int_{k_1\vartheta}^{(k_1+1)\vartheta} {|t-u|^{-\frac{1}{2}+\frac{\delta}{2}}|y-y_1|^{-2-\delta}}du
\\&{|\bar{t}_2-\bar{t}_1|^{\frac{\delta}{4}}}{(\bar{t}_2-\tilde{u})^{-\frac{1}{2}+\frac{\delta}{4}}
|\bar{y}_2-y_1|^{-2-\delta}}d\tilde{u}\frac{1}{\vartheta}1_{\{|y_1|\leq C\}}(\frac{\vartheta}{\varepsilon^2}\wedge1)dy_1
.\endaligned$$
Moreover, by Young's inequality we have $${|\bar{t}_2-\bar{t}_1|^{\frac{\delta}{4}}}{(\bar{t}_2-\tilde{u})^{-\frac{1}{2}+\frac{\delta}{4}}}{|t-u|^{-\frac{1}{2}+\frac{\delta}{2}}}\leq {|\bar{t}_2-\bar{t}_1|^{\frac{\delta}{2}}}{|t-u|^{-1+\delta}}+{(\bar{t}_2-\tilde{u})^{-1+\frac{\delta}{2}}},$$ which combined with [Hai14, Lemma 10.14] implies that
$$\aligned J\lesssim&|\bar{t}_2-\bar{t}_1|^{\frac{\delta}{2}}|y-\bar{y}_2|^{-1-2\delta}.\endaligned$$
Furthermore, we estimate $J$ in the following three cases:

Case I: $t-\bar{t}_2\geq 2C_0\varepsilon^2+4\vartheta$. Similar arguments as in  the proof of (4.6) imply that
$t-u\geq\frac{t-\bar{t}_2}{2}$, which combined with (4.7) and (4.13) implies that for $\delta>0$
$$\aligned J\lesssim&\sum_{k_1=[\frac{\bar{t}_1}{\vartheta}]+1}^{[\frac{\bar{t}_2}{\vartheta}]-1}\int\int_{k_1\vartheta}^{(k_1+1)\vartheta} {(t-\bar{t}_2)^{-\frac{1}{2}-\delta}|y-y_1|^{-2+2\delta}}
{(\bar{t}_2-\tilde{u})^{-1+\delta/2}|\bar{y}_2-y_1|^{-1-\delta}}d\tilde{u}1_{|y_1|\leq C}dy_1
\\\lesssim&|t-\bar{t}_2|^{-\frac{1}{2}-\delta}|\bar{t}_2-\bar{t}_1|^{\frac{\delta}{2}}.\endaligned$$

Case II: $\bar{t}_2-t\geq 2C_0\varepsilon^2+4\vartheta$.  Similarly as in Case I, we have $\bar{t}_2-\tilde{u}\geq\frac{\bar{t}_2-t}{2}$, which combined  with $|\bar{t}_2-\bar{t}_1|\geq |\bar{t}_2-\tilde{u}|$ implies that
$$\aligned J\lesssim&\sum_{k=-[\frac{C}{\vartheta}]}^{[\frac{t}{\vartheta}]-1}\int\int_{k\vartheta}^{(k+1)\vartheta} {(t-u)^{-1+\frac{\delta}{4}}|y-y_1|^{-1-\frac{\delta}{2}}}
{|\bar{t}_2-\bar{t}_1|^{\frac{\delta}{2}}}{(\bar{t}_2-t)^{-\frac{1}{2}-\delta}|\bar{y}_2-y_1|^{-2+\delta}}du1_{|y_1|\leq C}dy_1
\\\lesssim&|t-\bar{t}_2|^{-\frac{1}{2}-\delta}|\bar{t}_2-\bar{t}_1|^{\frac{\delta}{2}}.\endaligned$$

Case III: $|\bar{t}_2-t|\leq 2C_0\varepsilon^2+4\vartheta$. In this case
$$\frac{1}{\vartheta}\int_{k\vartheta}^{(k+1)\vartheta}{(t-u)^{-\frac{1}{2}-{\delta}}}du\lesssim \vartheta^{-\frac{1}{2}-{\delta}},$$
which combined with [Hai14, Lemma 10.14] implies that
$$\aligned J\lesssim&\sum_{k_1=[\frac{\bar{t}_1}{\vartheta}]+1}^{[\frac{\bar{t}_2}{\vartheta}]-1}\int\int_{k_1\vartheta}^{(k_1+1)\vartheta} {\vartheta^{-\frac{1}{2}-\delta}|y-y_1|^{-2+2\delta}}
{(\bar{t}_2-\tilde{u})^{-1+\delta/2}|\bar{y}_2-y_1|^{-1-\delta}}d\tilde{u}1_{|y_1|\leq C}dy_1\\\lesssim&\vartheta^{-\frac{1}{2}-\delta}|\bar{t}_2-\bar{t}_1|^{\frac{\delta}{2}}.\endaligned$$
On the other hand, a similar argument as in (4.11) and $|\bar{t}_2-\bar{t}_1|\geq |\bar{t}_2-\tilde{u}|$ imply that
$$\aligned J\lesssim&\sum_{k=-[\frac{C}{\vartheta}]}^{[\frac{t}{\vartheta}]-1}\sum_{k_1=[\frac{\bar{t}_1}{\vartheta}]+1, |k_1-k|\leq C_0\frac{\varepsilon^2}{\vartheta}+1}^{[\frac{\bar{t}_2}{\vartheta}]-1}\int\int_{k\vartheta}^{(k+1)\vartheta}\int_{k_1\vartheta}^{(k_1+1)\vartheta} {|t-u|^{-\frac{1}{2}+\frac{\delta}{8}}|y-y_1|^{-\frac{3}{2}+\frac{\delta}{4}}}\varepsilon^{-1-2\delta}du
\\&{|\bar{t}_2-\bar{t}_1|^{\delta/2}}{(\bar{t}_2-\tilde{u})^{-\frac{1}{2}+\frac{\delta}{8}}
|\bar{y}_2-y_1|^{-\frac{3}{2}+\frac{\delta}{4}}}d\tilde{u}\frac{1}{\vartheta}1_{\{|y_1|\leq C\}}(\frac{\vartheta}{\varepsilon^2}\wedge1)dy_1\\\lesssim&\varepsilon^{-1-2\delta}|\bar{t}_2-\bar{t}_1|^{\frac{\delta}{2}}.\endaligned$$
Combining the above estimates we obtain that
$$\aligned J\lesssim&\|\bar{z}_1-\bar{z}_2\|_\mathfrak{s}^\delta\|z-\bar{z}_2\|_\mathfrak{s}^{-1-2\delta}.\endaligned$$
 Thus (i) follows. Combining the arguments in (i) and the proof for Lemma 4.2 (iii), we can deduce (ii)  easily.
$\hfill\Box$

\vskip.10in

{\textbf{Appendix. Proof of Theorem 3.7}}
\vskip.10in

 By [Hai14, Theorem 10.7] we only need to show that the renormalised model converges for those elements $\tau\in\mathcal{F}_F$ with non-positive homogeneity. In the case of the dynamical $\Phi_3^4$ model, these elements are given by
$$\mathcal{F}_-=\{\Xi, \Psi,\Psi^2,\Psi^3,\Psi^2X_i, \mathcal{I}(\Psi^3)\Psi, \mathcal{I}(\Psi^2)\Psi^2, \mathcal{I}(\Psi^3)\Psi^2\}.$$

By [Hai14, Theorem 10.7] it is sufficient to prove that for $\tau\in\mathcal{F}_-$ with $|\tau|_\mathfrak{s}<0$,  any test function $\varphi\in\mathcal{B}_r$ and every $x\in \mathbb{R}^4$,
and  for some $0<\kappa<-\frac{5}{2}-\alpha$,
$$\mathbf{E}|(\hat{\Pi}_x^{(\varepsilon,\vartheta)}\tau)(\varphi_x^\lambda)|^2\lesssim \lambda^{2|\tau|_{\mathfrak{s}}+\kappa},\quad \mathbf{E}|(\hat{\Pi}_x\tau-\hat{\Pi}_x^{(\varepsilon,\vartheta)}\tau)(\varphi_x^\lambda)|^2\lesssim (\varepsilon^{2\kappa}+\vartheta^\kappa)\lambda^{2|\tau|_{\mathfrak{s}}+\kappa},\eqno(A.1)$$
where $\hat{\Pi}_x\tau$ is obtained as in the proof of  [Hai14, Theorem 10.22].
Since the map $\varphi\mapsto(\hat{\Pi}_x^{(\varepsilon,\vartheta)}\tau)(\varphi)$ is linear, we can find some functions $\hat{\mathcal{W}}_x^{(\varepsilon,\vartheta;k)}\tau$ with $(\hat{\mathcal{W}}_x^{(\varepsilon,\vartheta;k)}\tau)(y)\in L^2(\mathbb{R}\times \mathbb{T}^3)^{\otimes k}$, for $y\in\mathbb{R}^4$ and satisfying
$$(\hat{\Pi}_x^{(\varepsilon,\vartheta)}\tau)(\varphi)=\sum_{k\leq \|\tau\|}I_k\bigg(\int\varphi(y)(\hat{\mathcal{W}}_x^{(\varepsilon,\vartheta;k)}\tau)(y)dy\bigg),$$
where $\|\tau\|$ denotes the number of occurrences of $\Xi$ in the expression $\tau$ and $I_k$ is defined as in [Hai14, Section 10.1].
We also use the following notation as in [Hai14, Section 10]:
 $$(\hat{\Pi}_x\tau)(\varphi)=\sum_{k\leq \|\tau\|}I_k\bigg(\int\varphi(y)(\hat{\mathcal{W}}^{(k)}_x\tau)(y)dy\bigg),$$
 where $\hat{\mathcal{W}}^{(k)}_x\tau\in L^2(\mathbb{R}\times \mathbb{T}^3)^{\otimes k}$. By [Hai14, Proposition 10.11] we know that to obtain (A.1) it suffices to estimate the terms $|\langle(\hat{\mathcal{W}}_x^{(\varepsilon,\vartheta;k)}\tau)(z),(\hat{\mathcal{W}}_x^{(\varepsilon,\vartheta;k)}\tau)(\bar{z})\rangle|$ and $|\langle(\delta\hat{\mathcal{W}}_x^{(\varepsilon,\vartheta;k)}\tau)(z),(\delta\hat{\mathcal{W}}_x^{(\varepsilon,\vartheta;k)}\tau)(\bar{z})\rangle|$,
 where $\delta\hat{\mathcal{W}}_x^{(\varepsilon,\vartheta;k)}\tau=\hat{\mathcal{W}}^{(\varepsilon,\vartheta;k)}_x\tau-\hat{\mathcal{W}}_x^{(k)}\tau$.

 For $\tau=\mathcal{I}(\Xi)=\Psi$ we have
$$\aligned ({\hat{\Pi}}^{(\varepsilon,\vartheta)}\Psi)(z)=&K*\xi_{\varepsilon,\vartheta}(z)=\int K_{\varepsilon,\vartheta}(z,z_1)\xi(z_1)dz_1,\endaligned$$
which implies that
$$(\hat{\mathcal{W}}_x^{(\varepsilon,\vartheta;1)}\Psi)(z,z_1)=K_{\varepsilon,\vartheta}(z,z_1).$$
For $\varphi$ smooth and $x\in\mathbb{R}^4$ we have that
$$\aligned &E|\langle K*\xi_{\varepsilon,\vartheta},\varphi_{x}^\lambda\rangle|^{2}=\int\int f^{(\varepsilon,\vartheta)}(z,\bar{z})
\varphi_{x}^\lambda(z)\varphi_{x}^\lambda(\bar{z})dzd\bar{z}.\endaligned$$

In the following we use \includegraphics[height=0.5cm]{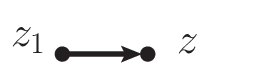} to represent a factor $K(z-z_1)$  and \includegraphics[height=0.5cm]{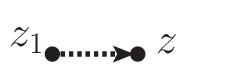} to represent  $K_{\varepsilon,\vartheta}(z,z_1)$. We also use the convention that if a vertex is drawn in grey, then the corresponding variable is integrated out.
Now we have
$$f^{(\varepsilon,\vartheta)}(z,\bar{z})=\includegraphics[height=1cm]{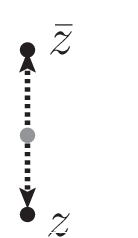}.$$
By Lemma 4.2 (i) we obtain that
$$\aligned &|\langle(\hat{\mathcal{W}}_x^{(\varepsilon,\vartheta;1)}\Psi)(z),(\hat{\mathcal{W}}_x^{(\varepsilon,\vartheta;1)}\Psi)(\bar{z})\rangle|
=|\includegraphics[height=1cm]{0013.eps}|\lesssim
\|z-\bar{z}\|_\mathfrak{s}^{-1-\delta},\endaligned$$
holds uniformly over $\varepsilon,\vartheta\in(0,1)$. Now for $\hat{\Pi}_x\Psi=K*\xi$ as in the proof of [Hai14, Theorem 10.22]  we also have
$$\aligned (\delta\hat{\mathcal{W}}_x^{(\varepsilon,\vartheta;1)}\Psi)(z,z_1)=&
K_{\varepsilon,\vartheta}(z,z_1)-K(z-z_1)\\=&
(K_{\varepsilon,\vartheta}^{(1)}(z,z_1)-K_\varepsilon^{(1)}(z,z_1))+(K_{\varepsilon,\vartheta}^{(2)}(z,z_1)-K_\varepsilon^{(2)}(z,z_1))
\\&+(K_{\varepsilon}(z-z_1)-K(z-z_1))\\:=&(\delta\hat{\mathcal{W}}_x^{(\varepsilon,\vartheta;11)}\Psi)(z,z_1)
+(\delta\hat{\mathcal{W}}_x^{(\varepsilon,\vartheta;12)}\Psi)(z,z_1)+(\delta\hat{\mathcal{W}}_x^{(\varepsilon,\vartheta;13)}\Psi)(z,z_1).\endaligned$$
By Lemma 4.2 (iii) and [Hai14, Lemmas 10.14, 10.17] we have that for $i=1,2$
$$|\langle(\delta\hat{\mathcal{W}}_x^{(\varepsilon,\vartheta;1i)}\Psi)(z),(\delta\hat{\mathcal{W}}_x^{(\varepsilon,\vartheta;1i)}\Psi)(\bar{z})\rangle|
\lesssim \vartheta^{\kappa}\|z-\bar{z}\|_\mathfrak{s}^{-1-\delta-2\kappa},\eqno(A.2)$$
and
$$|\langle(\delta\hat{\mathcal{W}}_x^{(\varepsilon,\vartheta;13)}\Psi)(z),(\delta\hat{\mathcal{W}}_x^{(\varepsilon,\vartheta;13)}\Psi)(\bar{z})\rangle|
\lesssim \varepsilon^{2\kappa}\|z-\bar{z}\|_\mathfrak{s}^{-1-\delta-2\kappa}, \eqno(A.3)$$
holds uniformly over $\varepsilon,\vartheta\in(0,1]$ provided that $0<\kappa<1$ and that $0<\delta<1$, from which we deduce (A.1) for $\tau=\Psi$ easily. In the following we use $|\includegraphics[height=0.5cm]{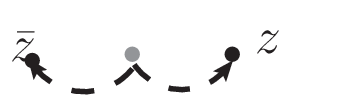}|$ to represent  $$\sum_{i=1}^3|\langle(\delta\hat{\mathcal{W}}_x^{(\varepsilon,\vartheta;1i)}\Psi)(z),(\delta\hat{\mathcal{W}}_x^{(\varepsilon,\vartheta;1i)}\Psi)(\bar{z})\rangle|.$$
By (A.2) and (A.3)  we have
$$|\includegraphics[height=1cm]{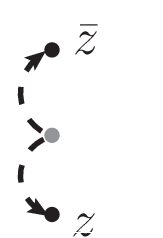}|\lesssim (\vartheta^\kappa+\varepsilon^{2\kappa})\|z-\bar{z}\|_\mathfrak{s}^{-1-\delta-2\kappa}.$$

For $\tau=\Psi^2$ we could choose for $z=(t,y)$
$$C_1^{(\varepsilon,\vartheta)}(t)=\int K_{\varepsilon,\vartheta}(z,z_1)^2dz_1.\eqno(A.4)$$
Here, since $\rho=\rho_1\rho_2$, we can easily deduce that $C_1^{(\varepsilon,\vartheta)}$ only depends on $t$.
We obtain that
$$\aligned &(\hat{\mathcal{W}}_x^{(\varepsilon,\vartheta,2)}\Psi^2)(z)=\includegraphics[height=0.7cm]{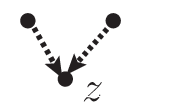}.\endaligned$$
By Lemma 4.2 we have that for every $\delta>0$
$$\aligned &|\langle(\hat{\mathcal{W}}_x^{(\varepsilon,\vartheta,2)}\Psi^2)(z),(\hat{\mathcal{W}}_x^{(\varepsilon,\vartheta,2)}\Psi^2)(\bar{z})\rangle|
= f^{(\varepsilon,\vartheta)}(z,\bar{z})^2\lesssim {\|z-\bar{z}\|_{\mathfrak{s}}^{-2-\delta}},\endaligned$$
holds uniformly over $\varepsilon,\vartheta\in(0,1]$. As in the proof of [Hai14, Theorem 10.22] $\hat{\mathcal{W}}_x^{(2)}\Psi^2(z; z_1, z_2) =
K(z,z_1)K(z, z_2)$.
By  (A.2), (A.3)  and Lemma 4.2 we have that
$$\aligned &|\langle(\delta\hat{\mathcal{W}}_x^{(\varepsilon,\vartheta,2)}\Psi^2)(z),(\delta\hat{\mathcal{W}}_x^{(\varepsilon,\vartheta,2)}\Psi^2)(\bar{z})\rangle|
\lesssim(\vartheta^\kappa+\varepsilon^{2\kappa})
{\|z-\bar{z}\|_\mathfrak{s}^{-2-2\kappa-{\delta}}},\endaligned$$
holds uniformly over $\varepsilon,\vartheta\in(0,1]$, provided that $0<\kappa<1$ and that $1>\delta>0$, which implies that (A.1) holds for $\tau=\Psi^2$.

Similar arguments also imply  that (A.1) holds for $\tau=\Psi^3$.

Regarding $\tau=\Psi^2X_i$ the corresponding bound follows from those for $\tau=\Psi^2$.

Now for $\tau=\mathcal{I}(\Psi^3)\Psi$ we have
$$(\hat{\Pi}_x^{(\varepsilon,\vartheta)}\tau)(z)=(\hat{\Pi}_x^{(\varepsilon,\vartheta)}\Psi)(z)[ K*(\hat{\Pi}_x^{(\varepsilon,\vartheta)}\Psi^3)(z)-K*(\hat{\Pi}_x^{(\varepsilon,\vartheta)}\Psi^3)(x)].$$
For the term in the fourth Wiener chaos we have
$$(\hat{\mathcal{W}}_x^{(\varepsilon,\vartheta,4)}\tau)(z)=\includegraphics[height=1cm]{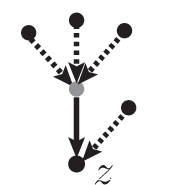}-\includegraphics[height=1cm]{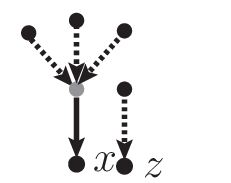}.$$
We have the following estimates:
$$\aligned &|\langle(\hat{\mathcal{W}}_x^{(\varepsilon,\vartheta,4)}\tau)(z),(\hat{\mathcal{W}}_x^{(\varepsilon,\vartheta,4)}\tau)(\bar{z})\rangle|
\\\lesssim&\int\int\|z-\bar{z}\|_\mathfrak{s}^{-1-\delta}
|K(z-z_1)-K(x-z_1)|
\|z_1-\bar{z}_1\|_\mathfrak{s}^{-3-\delta}\\&|K(\bar{z}-\bar{z}_1)-K(x-\bar{z}_1)|
dz_1d\bar{z}_1,\endaligned\eqno(A.5)$$
where we used Lemma 4.2 to obtain the estimate. We now use [Hai14, Lemma 10.18] to control $|K(z-z_1)-K(x-z_1)|$ by $\|z-x\|_\mathfrak{s}^{\frac{1}{2}-{\delta}} \bigg({\|z-z_1\|_{\mathfrak{s}}^{-3.5+\delta}}+{\|x-z_1\|_{\mathfrak{s}}^{-3.5+\delta}}\bigg)$ with $0<\delta<\frac{1}{2}$ and obtain that
$$\aligned &|\langle(\hat{\mathcal{W}}_x^{(\varepsilon,\vartheta,4)}\tau)(z),(\hat{\mathcal{W}}_x^{(\varepsilon,\vartheta,4)}\tau)(\bar{z})\rangle|\\
\lesssim&\|z-\bar{z}\|_\mathfrak{s}^{-1-\delta}
\|z-x\|_\mathfrak{s}^{\frac{1}{2}-{\delta}}\|\bar{z}-x\|_\mathfrak{s}^{\frac{1}{2}-{\delta}}(G(z-x)+G(\bar{z}-x)+G(z-\bar{z})+G(0))
,\endaligned\eqno(A.6)$$
holds uniformly over $\varepsilon,\vartheta\in(0,1]$, where the function $G$ is a bounded function  given by
$$G(z-\bar{z})=\includegraphics[height=0.7cm]{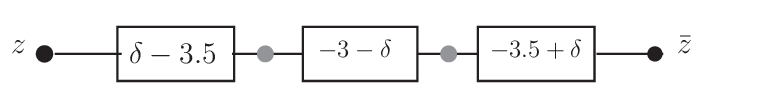}.\eqno(A.7)$$
  Here as in [Hai14, Theorem 10.22] we also use the notation $\includegraphics[height=0.5cm]{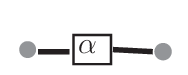}$ to represent $\|z-\bar{z}\|_\mathfrak{s}^{\alpha}1_{\|z-\bar{z}\|_\mathfrak{s}\leq C}$ for a constant $C$.
 Choose $\hat{\mathcal{W}}_x^{(4)}\tau$ as $\hat{\mathcal{W}}_x^{(\varepsilon,\vartheta,4)}\tau$ with  each instance of $K_{\varepsilon,\vartheta}$ replaced by $K$,
 which is the same as in the proof of [Hai14, Theorem 10.22]. By  (A.2), (A.3),
 Lemma 4.2 (i) and [Hai14, Lemma 10.18] we deduce that
$$\aligned |\langle(\delta\hat{\mathcal{W}}_x^{(\varepsilon,\vartheta,4)}\tau)(z),(\delta\hat{\mathcal{W}}_x^{(\varepsilon,\vartheta,4)}\tau)(\bar{z})\rangle|\lesssim&(\vartheta^{\kappa}
+\varepsilon^{2\kappa})[\|z-\bar{z}\|_\mathfrak{s}^{-1-\delta}
\|z-x\|_\mathfrak{s}^{\frac{1}{2}-{\kappa}-\delta}\|\bar{z}-x\|_\mathfrak{s}^{\frac{1}{2}-{\kappa}-\delta}
\\&+\|z-\bar{z}\|_\mathfrak{s}^{-1-{2\kappa}-\delta}
\|z-x\|_\mathfrak{s}^{\frac{1}{2}-\delta}\|\bar{z}-x\|_\mathfrak{s}^{\frac{1}{2}-\delta}],\endaligned$$
 holds uniformly over $\varepsilon,\vartheta\in(0,1]$, provided that $0<\kappa<1$ and that  $1>\delta>0$. For the term in the second Wiener chaos, we also have the following identity:
$$(\hat{\mathcal{W}}_x^{(\varepsilon,\vartheta,2)}\tau)(z)=3(\includegraphics[height=1cm]{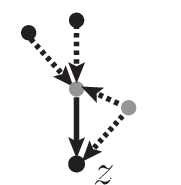}-\includegraphics[height=1cm]{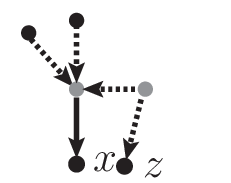})
:=3((\hat{\mathcal{W}}_x^{(\varepsilon,\vartheta,21)}\tau)(z)-(\hat{\mathcal{W}}_x^{(\varepsilon,\vartheta,22)}\tau)(z)).$$
For $\hat{\mathcal{W}}_x^{(\varepsilon,\vartheta,21)}\tau$  we have that for every $\delta>0$
$$\aligned &|\langle(\hat{\mathcal{W}}_x^{(\varepsilon,\vartheta,21)}\tau)(z),(\hat{\mathcal{W}}_x^{(\varepsilon,\vartheta,21)}\tau)(\bar{z})\rangle|\\
\lesssim&|\includegraphics[height=0.7cm]{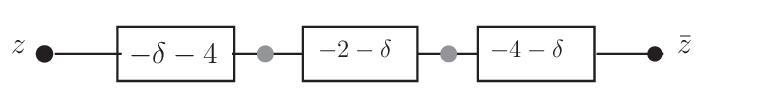}|
\lesssim\|z-\bar{z}\|_\mathfrak{s}^{-{3\delta}}
,\endaligned$$
where we used Lemma 4.2 (i) in the first inequality and [Hai14, Lemma 10.14] in the last inequality. Choose $\hat{\mathcal{W}}_x^{(21)}\tau$ as $\hat{\mathcal{W}}_x^{(\varepsilon,\vartheta,21)}\tau$ with each instance of $K_{\varepsilon,\vartheta}$ replaced by $K$,
 which is the same as in the proof of [Hai14, Theorem 10.22]. By Lemmas 4.2 and (A.2), (A.3)  we have that
$$\aligned &|\langle(\delta\hat{\mathcal{W}}_x^{(\varepsilon,\vartheta,21)}\tau)(z),(\delta\hat{\mathcal{W}}_x^{(\varepsilon,\vartheta,21)}\tau)(\bar{z})\rangle|
\lesssim(\vartheta^{\kappa}+\varepsilon^{2\kappa})\|z-\bar{z}\|_\mathfrak{s}^{-{2\kappa}-3\delta}
,\endaligned$$
holds uniformly over $\varepsilon,\vartheta\in(0,1]$, provided that $0<\kappa<1$ and that  $0<\delta<1$.
For $\hat{\mathcal{W}}_x^{(\varepsilon,\vartheta,22)}\tau$ by Lemma 4.2 we have that
$$\aligned &|\langle(\hat{\mathcal{W}}_x^{(\varepsilon,\vartheta,22)}\tau)(z),(\hat{\mathcal{W}}_x^{(\varepsilon,\vartheta,22)}\tau)(\bar{z})\rangle|\\
\lesssim& |\includegraphics[height=1.5cm]{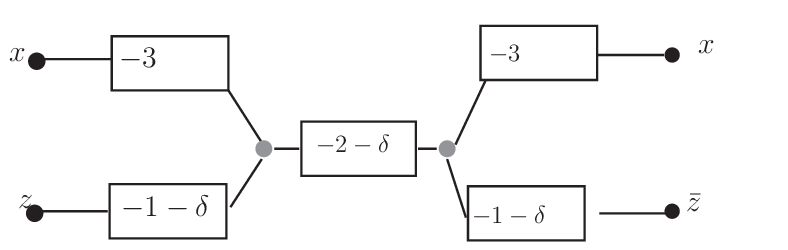}|
\\\lesssim&\|z-x\|_\mathfrak{s}^{-{4\delta}}(G^1(z-x)+G^1(\bar{z}-x)+G^1(z-\bar{z})+G^1(0))
,\endaligned$$
holds uniformly over $\varepsilon,\vartheta\in(0,1]$, provided that $0<\kappa<1$ and that  $0<\delta<1$, where we used [Hai14, (10.37)] in the last inequality and that the function $G^1$ is a bounded function  given by
 $$G^1(z-\bar{z})=\includegraphics[height=0.7cm]{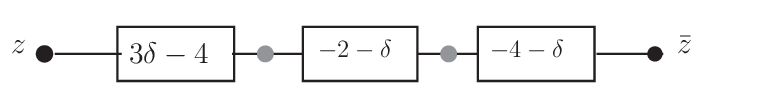}.$$
 Define $\hat{\mathcal{W}}_x^{(22)}\tau$ as $\hat{\mathcal{W}}_x^{(\varepsilon,\vartheta,22)}\tau$ with each instance of $K_{\varepsilon,\vartheta}$ replaced by $K$. For the difference by Lemma 4.2 and (A.2), (A.3) we have that
$$\aligned &|\langle(\delta\hat{\mathcal{W}}_x^{(\varepsilon,\vartheta,22)}\tau)(z),(\delta\hat{\mathcal{W}}_x^{(\varepsilon,\vartheta,22)}\tau)(\bar{z})\rangle|
\lesssim(\vartheta^{\kappa}+\varepsilon^{2\kappa})\|z-x\|_\mathfrak{s}^{-{2\kappa}-4\delta}
,\endaligned$$
holds uniformly over $\varepsilon,\vartheta\in(0,1]$, provided that $0<\kappa<1$ and that  $0<\delta<1$.
Combining all the estimates above we obtain that (A.1)  holds for $\tau=\mathcal{I}(\Psi^3)\Psi$.

Now we come to the case $\tau=\mathcal{I}(\Psi^2)\Psi^2$. We have for $z=(t,y)$
$$(\hat{\Pi}_x^{(\varepsilon,\vartheta)}\tau)(z)=(\hat{\Pi}_x^{(\varepsilon,\vartheta)}\Psi^2)(z)[ K*(\hat{\Pi}_x^{(\varepsilon,\vartheta)}\Psi^2)(z)-K*(\hat{\Pi}_x^{(\varepsilon,\vartheta)}\Psi^2)(x)]-C_2^{(\varepsilon,\vartheta)}(t).$$
For the term in the fourth Wiener chaos, we have
$$(\hat{\mathcal{W}}_x^{(\varepsilon,\vartheta,4)}\tau)(z)=\includegraphics[height=1cm]{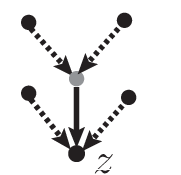}-\includegraphics[height=1cm]{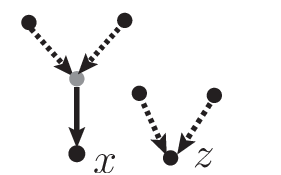}.$$
By similar calculations as in (A.5), (A.6) we have that
$$\aligned &|\langle(\hat{\mathcal{W}}_x^{(\varepsilon,\vartheta,4)}\tau)(z),(\hat{\mathcal{W}}_x^{(\varepsilon,\vartheta,4)}\tau)(\bar{z})\rangle|\\\lesssim&\|z-\bar{z}\|_\mathfrak{s}^{-2-\delta}\|z-x\|_\mathfrak{s}^{1-{\delta}}\|\bar{z}-x\|_\mathfrak{s}^{1-{\delta}}(G^2(z-x)+G^2(\bar{z}-x)+G^2(z-\bar{z})+G^2(0))\endaligned$$
holds uniformly over $\varepsilon,\vartheta\in(0,1]$, where we  apply [Hai14, Lemma 10.18] to control $|K(z-z_1)-K(x-z_1)|$ by $\|z-x\|_\mathfrak{s}^{1-{\delta}} \bigg({\|z-z_1\|_{\mathfrak{s}}^{-4+\delta}}$ $+{\|x-z_1\|_{\mathfrak{s}}^{-4+\delta}}\bigg)$ for $0<\delta<1$ and the function $G^2$ is a bounded function  given by
$$G^2(z-\bar{z})=\includegraphics[height=0.7cm]{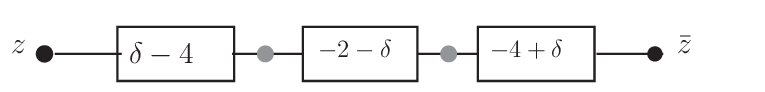}.$$
Choose $\hat{\mathcal{W}}_x^{(4)}\tau$ as $\hat{\mathcal{W}}_x^{(\varepsilon,\vartheta,4)}\tau$ with each instance of $K_{\varepsilon,\vartheta}$ replaced by $K$, which is the same as in the proof of [Hai14, Theorem 10.22].
Similarly, by Lemma 4.2 and (A.2), (A.3) we have that
$$\aligned &|\langle(\delta\hat{\mathcal{W}}_x^{(\varepsilon,\vartheta,4)}\tau)(z),(\delta\hat{\mathcal{W}}_x^{(\varepsilon,\vartheta,4)}\tau)(\bar{z})\rangle|\\
\lesssim&(\vartheta^{\kappa}+\varepsilon^{2\kappa})
[\|z-\bar{z}\|_\mathfrak{s}^{-2-\delta}\|z-x\|_\mathfrak{s}^{1-{\kappa}-\delta}\|\bar{z}-x\|_\mathfrak{s}^{1-{\kappa}-\delta}
+\|z-\bar{z}\|_\mathfrak{s}^{-2-{2\kappa}-\delta}
\|z-x\|_\mathfrak{s}^{1-\delta}\|\bar{z}-x\|_\mathfrak{s}^{1-\delta}],\endaligned$$
holds uniformly over $\varepsilon,\vartheta\in(0,1]$, provided that $0<\kappa<1$ and that  $1>\delta>0$. For the term in the second Wiener chaos, we have the following identity
$$\aligned(\hat{\mathcal{W}}_x^{(\varepsilon,\vartheta,2)}\tau)(z)=4(\includegraphics[height=1cm]{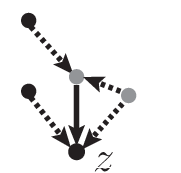}-\includegraphics[height=1cm]{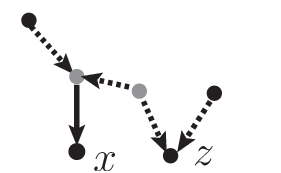}).\endaligned$$
Then by Lemma 4.2 we obtain that
$$\aligned &|\langle(\hat{\mathcal{W}}_x^{(\varepsilon,\vartheta,2)}\tau)(z),(\hat{\mathcal{W}}_x^{(\varepsilon,\vartheta,2)}\tau)(\bar{z})\rangle|\\
\lesssim&\int\int\|z-\bar{z}\|_\mathfrak{s}^{-1-\delta}{\|z-z_1\|_\mathfrak{s}^{-1-\delta}}|K(z-z_1)-K(x-z_1)|
{\|z_1-\bar{z}_1\|_\mathfrak{s}^{-1-\delta}}\\&|K(\bar{z}-\bar{z}_1)-K(x-\bar{z}_1)|{\|\bar{z}-\bar{z}_1\|^{-1-\delta}_\mathfrak{s}}
dz_1d\bar{z}_1
\\\lesssim&\|z-\bar{z}\|_\mathfrak{s}^{-1-\delta}\|z-x\|_\mathfrak{s}^{\frac{1}{2}-2\delta}
\|\bar{z}-x\|_\mathfrak{s}^{\frac{1}{2}-2\delta}(G^3(z,\bar{z})+G^3(z,x)+G^3(x,\bar{z})+G^3(x,x)),\endaligned$$
holds uniformly over $\varepsilon,\vartheta\in(0,1]$. Here the function $G^3$ is a bounded function given by
$$G^3(a,b)=\includegraphics[height=1.5cm]{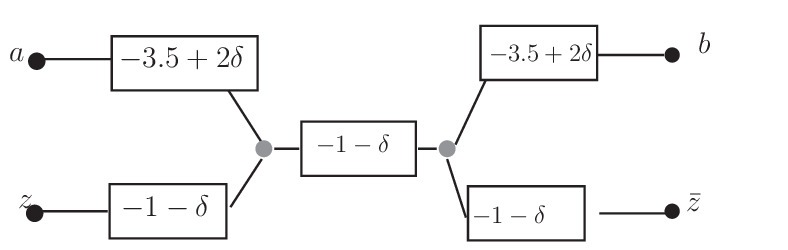},$$
and we used Young's inequality to obtain that $G^3$ is bounded. Choose $\hat{\mathcal{W}}_x^{(2)}\tau$ as $\hat{\mathcal{W}}_x^{(\varepsilon,\vartheta,2)}\tau$ with each instance of $K_{\varepsilon,\vartheta}$ replaced by $K$,
 which is the same as in the proof of [Hai14, Theorem 10.22].
Similarly, by Lemma 4.2 and (A.2), (A.3) we have that
$$\aligned &|\langle(\delta\hat{\mathcal{W}}_x^{(\varepsilon,\vartheta,2)}\tau)(z),(\delta\hat{\mathcal{W}}_x^{(\varepsilon,\vartheta,2)}\tau)(\bar{z})\rangle|\\\lesssim&
(\vartheta^{\kappa}+\varepsilon^{2\kappa})[\|z-\bar{z}\|_\mathfrak{s}^{-1-\delta}\|z-x\|_\mathfrak{s}^{\frac{1}{2}-2\delta-\kappa}
\|\bar{z}-x\|_\mathfrak{s}^{\frac{1}{2}-2\delta-\kappa}+\|z-\bar{z}\|_\mathfrak{s}^{-1-2\kappa-\delta}\|z-x\|_\mathfrak{s}^{\frac{1}{2}-2\delta}
\|\bar{z}-x\|_\mathfrak{s}^{\frac{1}{2}-2\delta}]
,\endaligned$$
holds uniformly over $\varepsilon,\vartheta\in(0,1]$, provided that $0<\kappa<1$ and that  $1>\delta>0$.
 We now turn to the component in the 0th Wiener chaos. For $z=(t,y)$, choose
$$C_2^{(\varepsilon,\vartheta)}(t)=2\int f^{(\varepsilon,\vartheta)}(z,z_1)^2K(z-z_1)dz_1=2\includegraphics[height=0.5cm]{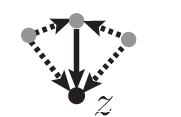}.\eqno(A.8)$$
Here $C_2^{(\varepsilon,\vartheta)}$  only depends on $t$.  We have
$$\aligned(\hat{\mathcal{W}}_x^{(\varepsilon,\vartheta,0)}\tau)(z)=-2\int f^{(\varepsilon,\vartheta)}(z,z_1)^2K(x-z_1)dz_1,\endaligned$$
which combined with Lemma 4.2 impies that
$$\aligned|(\hat{\mathcal{W}}_x^{(\varepsilon,\vartheta,0)}\tau)(z)|\lesssim&\int \|z-z_1\|_\mathfrak{s}^{-2-\delta}|K(x-z_1)|dz_1\lesssim \|z-x\|_\mathfrak{s}^{-\delta},\endaligned$$
for every $\delta>0$. Choose $\hat{\mathcal{W}}_x^{(0)}\tau$ as above with each instance of $K_{\varepsilon,\vartheta}$ replaced by $K$,
 which is the same as in the proof of [Hai14, Theorem 10.22].
Then Lemma 4.2 yields that
$$\aligned&|(\hat{\mathcal{W}}_x^{(\varepsilon,\vartheta,0)}\tau)(z)-(\hat{\mathcal{W}}_x^{(0)}\tau)(z)|\\\lesssim&(\varepsilon^{2\kappa}+\vartheta^{\kappa})\int \|z-z_1\|_\mathfrak{s}^{-2-\delta-2\kappa}|K(x-z_1)|dz_1\\\lesssim &(\varepsilon^{2\kappa}+\vartheta^{\kappa})\|z-x\|_\mathfrak{s}^{-\delta-2\kappa},\endaligned$$
holds uniformly over $\varepsilon,\vartheta\in(0,1]$, provided that $0<\kappa<1$ and that  $1>\delta>0$.

For $\tau=\mathcal{I}(\Psi^3)\Psi^2$, we have the following identity for $z=(t,y)$
$$(\hat{\Pi}_x^{(\varepsilon,\vartheta)}\tau)(z)=(\hat{\Pi}_x^{(\varepsilon,\vartheta)}\Psi^2)(z)[ K*(\hat{\Pi}_x^{(\varepsilon,\vartheta)}\Psi^3)(z)-K*(\hat{\Pi}_x^{(\varepsilon,\vartheta)}\Psi^3)(x)]-3C_2^{(\varepsilon,\vartheta)}(t)(\hat{\Pi}_x^{(\varepsilon,\vartheta)}\Psi)(z).$$
For the term in the fifth Wiener chaos, we have
$$(\hat{\mathcal{W}}_x^{(\varepsilon,\vartheta,5)}\tau)(z)=\includegraphics[height=1cm]{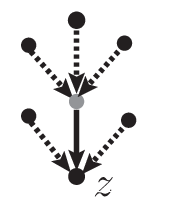}-\includegraphics[height=1cm]{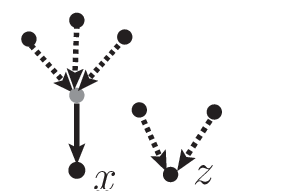}.$$
By  similar calculations as in (A.5) and (A.6) we obtain that
$$\aligned &|\langle(\hat{\mathcal{W}}_x^{(\varepsilon,\vartheta,5)}\tau)(z),(\hat{\mathcal{W}}_x^{(\varepsilon,\vartheta,5)}\tau)(\bar{z})\rangle|
\\\lesssim&\|z-\bar{z}\|_\mathfrak{s}^{-2-\delta}
\|z-x\|_\mathfrak{s}^{\frac{1}{2}-{\delta}}\|\bar{z}-x\|_\mathfrak{s}^{\frac{1}{2}-\delta}(G(z-x)+G(\bar{z}-x)+G(0)+G(z-\bar{z})),\endaligned$$
where  the function $G$ is given by (A.7).
Choose $\hat{\mathcal{W}}_x^{(5)}\tau$ as $\hat{\mathcal{W}}_x^{(\varepsilon,\vartheta,5)}\tau$ with each instance of $K_{\varepsilon,\vartheta}$ replaced by $K$,
 which is the same as in the proof of [Hai14, Theorem 10.22].
For the difference by Lemma 4.2 and (A.2), (A.3) we have similar estimates:
$$\aligned &|\langle(\delta\hat{\mathcal{W}}_x^{(\varepsilon,\vartheta,5)}\tau)(z),(\delta\hat{\mathcal{W}}_x^{(\varepsilon,\vartheta,5)}\tau)(\bar{z})\rangle|
\\\lesssim&(\vartheta^{\kappa}+\varepsilon^{2\kappa})[\|z-\bar{z}\|_\mathfrak{s}^{-2-\delta}
\|z-x\|_\mathfrak{s}^{\frac{1}{2}-{\kappa}-\delta}\|\bar{z}-x\|_\mathfrak{s}^{\frac{1}{2}-{\kappa}-\delta}
+\|z-\bar{z}\|_\mathfrak{s}^{-2-{2\kappa}-\delta}
\|z-x\|_\mathfrak{s}^{\frac{1}{2}-\delta}\|\bar{z}-x\|_\mathfrak{s}^{\frac{1}{2}-\delta}],\endaligned$$
which is valid uniformly over $\varepsilon,\vartheta\in(0,1]$, provided that $0<\kappa<1$ and that  $1>\delta>0$.

The component in the third Wiener chaos is very similar to what was obtained previously. Indeed, we have
$$(\hat{\mathcal{W}}_x^{(\varepsilon,\vartheta,3)}\tau)(z)=6(\includegraphics[height=1cm]{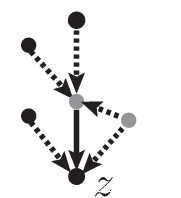}-\includegraphics[height=1cm]{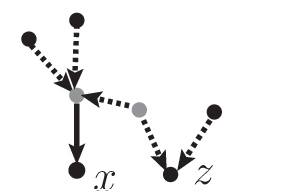})
:=6((\hat{\mathcal{W}}_x^{(\varepsilon,\vartheta,31)}\tau)(z)-(\hat{\mathcal{W}}_x^{(\varepsilon,\vartheta,32)}\tau)(z)).$$
Then we obtain that for every $\delta>0$
$$\aligned &|\langle(\hat{\mathcal{W}}_x^{(\varepsilon,\vartheta,31)}\tau)(z),(\hat{\mathcal{W}}_x^{(\varepsilon,\vartheta,31)}\tau)(\bar{z})\rangle|\\
\lesssim&\|z-\bar{z}\|_\mathfrak{s}^{-1-\delta}|\includegraphics[height=0.7cm]{0002.eps}|
\lesssim\|z-\bar{z}\|_\mathfrak{s}^{-1-4\delta},\endaligned$$
where we used Lemma 4.2 in the first inequality and [Hai14, Lemma 10.14] in the last inequality. Choose $\hat{\mathcal{W}}_x^{(31)}\tau$ as $\hat{\mathcal{W}}_x^{(\varepsilon,\vartheta,31)}\tau$ with each instance of $K_{\varepsilon,\vartheta}$ replaced by $K$,
 which is the same as in the proof of [Hai14, Theorem 10.22]. Similarly, by Lemma 4.2 and (A.2), (A.3) we have that
$$\aligned &|\langle(\delta\hat{\mathcal{W}}_x^{(\varepsilon,\vartheta,31)}\tau)(z),(\delta\hat{\mathcal{W}}_x^{(\varepsilon,\vartheta,31)}\tau)(\bar{z})\rangle|
\lesssim(\vartheta^{\kappa}+\varepsilon^{2\kappa})\|z-\bar{z}\|_\mathfrak{s}^{-1-{2\kappa}-4\delta},\endaligned$$
holds uniformly over $\varepsilon,\vartheta\in(0,1]$, provided that $0<\kappa<1$ and that  $1>\delta>0$.
Similarly, we obtain
$$\aligned &|\langle(\hat{\mathcal{W}}_x^{(\varepsilon,\vartheta,32)}\tau)(z),(\hat{\mathcal{W}}_x^{(\varepsilon,\vartheta,32)}\tau)(\bar{z})\rangle|
\\\lesssim&|\includegraphics[height=1.5cm]{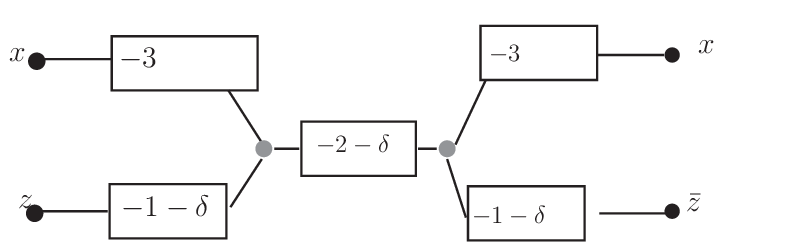}|\|z-\bar{z}\|_\mathfrak{s}^{-1-\delta}
\\\lesssim&[|\includegraphics[height=1.5cm]{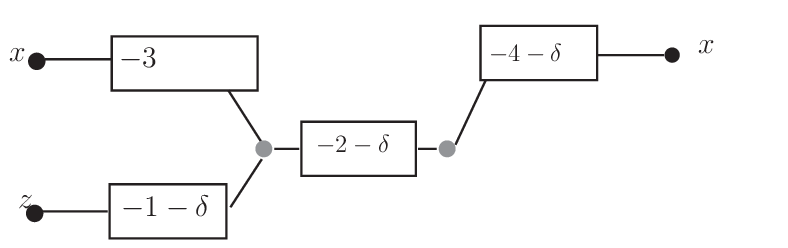}|+|\includegraphics[height=1.5cm]{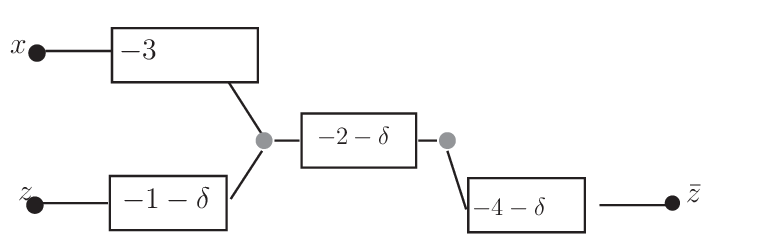}|]\|z-\bar{z}\|_\mathfrak{s}^{-1-\delta}
\\\lesssim&\|z-\bar{z}\|_\mathfrak{s}^{-1-\delta}(\|z-x\|_\mathfrak{s}^{-3\delta}
+\|z-\bar{z}\|_\mathfrak{s}^{-3\delta}),\endaligned$$
where we used Young's inequality in the second and last inequalities as well as [Hai14, Lemma 10.14] in the last inequality. By Lemma 4.2, (A.2), (A.3) and similar arguments as above we have
$$\aligned &|\langle(\delta\hat{\mathcal{W}}_x^{(\varepsilon,\vartheta,32)}\tau)(z),(\delta\hat{\mathcal{W}}_x^{(\varepsilon,\vartheta,32)}\tau)(\bar{z})\rangle|
\\\lesssim&(\vartheta^{\kappa}+\varepsilon^{2\kappa})[\|z-\bar{z}\|_\mathfrak{s}^{-1-\delta}(\|z-x\|_\mathfrak{s}^{-{2\kappa}-3\delta}
+\|z-\bar{z}\|_\mathfrak{s}^{-{2\kappa}-3\delta})+\|z-\bar{z}\|_\mathfrak{s}^{-1-2\kappa-\delta}(\|z-x\|_\mathfrak{s}^{-3\delta}
+\|z-\bar{z}\|_\mathfrak{s}^{-3\delta})],\endaligned$$
which is valid uniformly over $\varepsilon,\vartheta\in(0,1]$, provided that $0<\kappa<1$ and that  $1>\delta>0$.

We turn to the first Wiener chaos:
$$\aligned(\hat{\mathcal{W}}_x^{(\varepsilon,\vartheta,1)}\tau)(z)=&[(6\includegraphics[height=1cm]{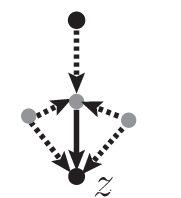}
-3C_2^{(\varepsilon,\vartheta)}(t)\includegraphics[height=0.5cm]{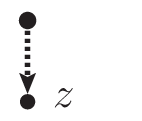})-6\includegraphics[height=1cm]{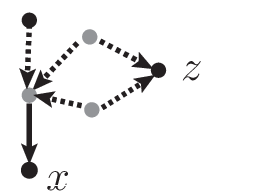}]
\\:=&
6[(\hat{\mathcal{W}}_x^{(\varepsilon,\vartheta,11)}\tau)(z)-(\hat{\mathcal{W}}_x^{(\varepsilon,\vartheta,12)}\tau)(z)].\endaligned$$
By Lemmas 4.2 and 4.3 we have that for every $\delta>0$
$$\aligned &|\langle(\hat{\mathcal{W}}_x^{(\varepsilon,\vartheta,11)}\tau)(z),(\hat{\mathcal{W}}_x^{(\varepsilon,\vartheta,11)}\tau)(\bar{z})\rangle|\\
\lesssim&\int\int |f^{(\varepsilon,\vartheta)}(z,z_1)^2K(z-z_1)[f^{(\varepsilon,\vartheta)}(z_1,\bar{z}_1)-f^{(\varepsilon,\vartheta)}(z,\bar{z}_1)
-f^{(\varepsilon,\vartheta)}(z_1,\bar{z})+f^{(\varepsilon,\vartheta)}(z,\bar{z})]\\&f^{(\varepsilon,\vartheta)}(\bar{z},\bar{z}_1)^2K(\bar{z}-\bar{z}_1)|
dz_1d\bar{z}_1\\\lesssim&\int\int\|z-z_1\|_\mathfrak{s}^{-5+\delta}
\|\bar{z}-\bar{z}_1\|_\mathfrak{s}^{-5+\delta}[\|z_1-\bar{z}_1\|^{-1-8\delta}_\mathfrak{s}+\|z-\bar{z}_1\|^{-1-8\delta}_\mathfrak{s}
+\|z_1-\bar{z}\|^{-1-8\delta}_\mathfrak{s}+\|z-\bar{z}\|_\mathfrak{s}^{-1-8\delta}]\\&1_{\{\|z_1\|_\mathfrak{s}\leq C\}}1_{\{\|\bar{z}_1\|_\mathfrak{s}\leq C\}}dz_1d\bar{z}_1
\\\lesssim&\|z-\bar{z}\|_\mathfrak{s}^{-1-8\delta}
,\endaligned\eqno(A.9)$$
where we used  interpolation in the second inequality and [Hai14, Lemma 10.14] in the last inequality.
Choose $\hat{\mathcal{W}}_x^{(11)}\tau(z,z_1)=\int (L(z-z_2)(K(z_2-z_1)-K(z-z_1))dz_2$ as in the proof of [Hai14, Theorem 10.22], where $L=(K*K)^2K$.

 Moreover, by  Lemma 4.3, interpolation and [Hai14, Lemmas 10.14, 10.17] we have
 $$\aligned&\sum_{i=1}^2|\langle (K^{(i)}_{\varepsilon,\vartheta}-K_\varepsilon^{(i)})(z,\cdot)-(K^{(i)}_{\varepsilon,\vartheta}-K_\varepsilon^{(i)})(z_1,\cdot),(K_{\varepsilon,\vartheta}^{(i)}-K^{(i)}_\varepsilon)(\bar{z},\cdot)
 -(K^{(i)}_{\varepsilon,\vartheta}-K_\varepsilon^{(i)})(\bar{z}_1,\cdot)\rangle|\\&+|\langle (K_{\varepsilon}-K_\varepsilon)(z,\cdot)-(K_{\varepsilon}-K_\varepsilon)(z_1,\cdot),(K_{\varepsilon}-K)(\bar{z},\cdot)
 -(K_{\varepsilon}-K)(\bar{z}_1,\cdot)\rangle|\\\lesssim& ( \varepsilon^{2\kappa}+\vartheta^\kappa)\|z-z_1\|_\mathfrak{s}^\delta\|\bar{z}-\bar{z}_1\|_\mathfrak{s}^\delta\\&(\|\bar{z}_1-z\|_\mathfrak{s}^{-1-4\delta-2\kappa}
 +\|\bar{z}-z\|_\mathfrak{s}^{-1-4\delta-2\kappa}+\|\bar{z}_1-z_1\|_\mathfrak{s}^{-1-4\delta-2\kappa}
 +\|\bar{z}-z_1\|_\mathfrak{s}^{-1-4\delta-2\kappa}),\endaligned$$
 which combined with similar arguments as those for (A.9) implies  the desired estimates for
$\delta\hat{\mathcal{W}}_x^{(\varepsilon,\vartheta,11)}\tau$.
Also by Lemma 4.2 we obtain that
$$\aligned &|\langle(\hat{\mathcal{W}}_x^{(\varepsilon,\vartheta,12)}\tau)(z),(\hat{\mathcal{W}}_x^{(\varepsilon,\vartheta,12)}\tau)(\bar{z})\rangle|\\
\lesssim&|\includegraphics[height=1.5cm]{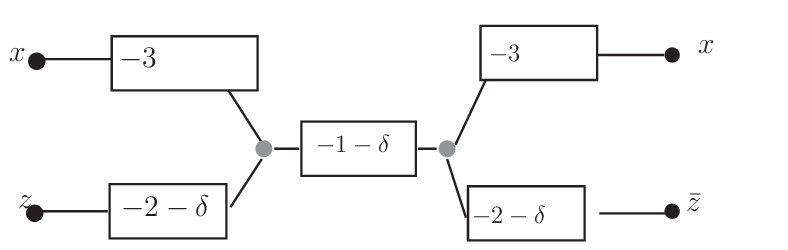}|
\\\lesssim&\|z-x\|_\mathfrak{s}^{-\frac{1}{2}-2\delta}\|\bar{z}-x\|_\mathfrak{s}^{-\frac{1}{2}-2\delta}(G^4(z-\bar{z})+G^4(z-x)+G^4(\bar{z}-x)+G^4(0))
,\endaligned$$
holds uniformly over $\varepsilon,\vartheta\in(0,1]$, provided that $0<\kappa<1$ and that  $1>\delta>0$. Here we used [Hai14, (10.37)] in the last inequality and that the function $G^4$ is a bounded function  given by
$$G^4(z-\bar{z})=\includegraphics[height=1cm]{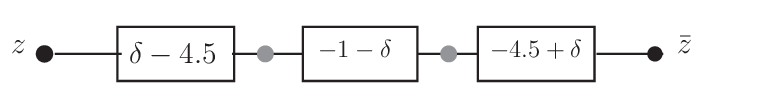}.$$

Similarly, by Lemma 4.2 and (A.2), (A.3), we have that
$$\aligned &|\langle(\delta\hat{\mathcal{W}}_x^{(\varepsilon,\vartheta,12)}\tau)(z),(\delta\hat{\mathcal{W}}_x^{(\varepsilon,\vartheta,12)}\tau)(\bar{z})\rangle|
\lesssim(\varepsilon^{2\kappa}+\vartheta^{\kappa})\|z-x\|_\mathfrak{s}^{\frac{-1-4\delta}{2}-\kappa}\|\bar{z}-x\|_\mathfrak{s}^{\frac{-1-4\delta}{2}-\kappa}
,\endaligned$$
holds uniformly over $\varepsilon,\vartheta\in(0,1]$, provided that $0<\kappa<1$ and that  $1>\delta>0$. Hence we conclude that (A.1)  holds for all $\tau\in\mathcal{F}^-$, which implies the results.

\end{document}